\documentclass[11pt]{article}

\usepackage{a4wide}
\DeclareMathAlphabet\mathbb{U}{fplmbb}{m}{n}
\usepackage{color}
\usepackage{url}
\usepackage{amssymb}%
\usepackage{pstricks,pst-node,pst-text,pst-3d}%
\usepackage{graphicx}
\usepackage{theorem}
\usepackage{psfrag}
\usepackage{floatflt}

\graphicspath{{Figures/}}
\input{psfig.sty}

\newtheorem{theorem}{Theorem}
{\theorembodyfont{\rmfamily}}
\newtheorem{lemma}{Lemma}
\newtheorem{proposition}{Proposition}
\newtheorem{corollary}{Corollary}
\newcommand{\qed}{\mbox{}\hspace*{\fill}{$\Box$}\medskip}

\def\MT{\widetilde{M}}     
\def\aa{\alpha}

\def\NN{\mathbb{N}}

\def\PP{\mathbb{P}}
\def\EE{\mathbb{E}}
\def\Eins{\mathbf{1}}
\def\MM{\mathcal{M}}
\def\MS{M^\sigma}          
\def\Int{{\sf Int}} 
\def\term#1{{\em #1}\marginpar{\raggedright{\small\it #1}}}
\def\ni{\noindent}
\def\MSD{M^{\lower1pt\hbox{\kern-.5pt$\scriptstyle \sigma$}\kern-1.5pt
            \raise2.5pt\hbox{$\scriptstyle *$}}}

\def\ITEMMACRO #1 ??? #2 ???{\par\vskip4pt\noindent%
\hangindent=#2em\setbox0\hbox{#1\kern4pt}%
\ifdim\wd0<\hangindent\setbox0\hbox to\hangindent{\hss#1\kern7pt}\fi%
\box0\ignorespaces}

\def\Item(#1){\ITEMMACRO {\rm (#1)} ??? 1.8 ???}

\def\BrackItem[#1]{\ITEMMACRO [#1] ??? 1.8 ???}
\def\proof{\par{\ni \sl Proof. }}

\def\bin#1#2{{#1\choose #2}}
\def\pbox#1#2#3{\parbox[b][#1cm][c]{#2cm}{#3}}
\def\cent#1#2{{\begin{center}\vspace{#1cm} #2\end{center}}}

\title{On the Number of Planar Orientations with Prescribed Degrees
\footnote{The conference version of this paper is to appear in the
proceedings of WG'07 (LNCS) under the title
"On the Number of $\alpha$-Orientations".} }
\author{Stefan Felsner \hskip2cm Florian Zickfeld\\[2mm]
\small Technische Universit\"at Berlin, Fachbereich Mathematik \\
\small Stra{\ss}e des 17. Juni 136, 10623 Berlin, Germany \\
\small \url{{felsner,zickfeld}@math.tu-berlin.de}}
\date{}

\begin{document}
\maketitle

\begin{abstract}
We deal with the asymptotic enumeration of combinatorial structures on planar
maps.  Prominent instances of such problems are the enumeration of
spanning trees, bipartite perfect matchings, and ice models. The notion
of orientations with out-degrees prescribed by a function
$\aa:V\rightarrow \NN$  unifies many different combinatorial
structures, including the afore mentioned.  We call these orientations
$\aa$-orientations. The main focus of this paper are bounds for the maximum
number of $\aa$-orientations that a planar map with $n$ vertices can
have, for different instances of $\aa$. We give examples of
triangulations with $2.37^n$ Schnyder woods, 
3-connected planar maps with $3.209^n$ Schnyder woods and inner
triangulations with $2.91^n$ bipolar orientations. These lower bounds
are accompanied by upper bounds of $3.56^n$, $8^n$ and $3.97^n$
respectively. We also show that for any planar map $M$ and any
$\alpha$ the number of $\alpha$-orientations is bounded from above by
$3.73^n$ and describe a family of maps which have at least $2.598^n$
$\alpha$-orientations.\bigskip

\ni \tt{AMS Math Subject Classification: 05A16, 05C20, 05C30}

\end{abstract}

\section{Introduction}
%

A \emph{planar map} is a planar graph together with a crossing-free
drawing in the plane.
Many different structures on planar maps have attracted the attention
of researchers. Among them are spanning trees, bipartite perfect
matchings (or more generally bipartite $f$-factors), Eulerian
orientations, Schnyder woods, bipolar orientations and 2-orientations
of quadrangulations.  The concept of orientations with prescribed
out-degrees is a quite general one. Remarkably, all the above
structures can be encoded as orientations with prescribed
out-degrees. Let a planar map $M$ with vertex set $V$ and a function
$\alpha: V \rightarrow \NN$ be given. An orientation $X$ of the edges
of $M$ is an $\alpha$-orientation if every vertex $v$ has
out-degree~$\alpha(v)$. For the sake of brevity, we refer to orientations
with prescribed out-degrees simply as $\aa$-orientations in this paper.

For some of the above mentioned structures it is not obvious how to
encode them as $\alpha$-orientations. For Schnyder woods on
triangulations the encoding by 3-orientations goes back to de~Fraysseix
and de~Mendez~\cite{fraysseixmendez01}. For bipolar orientations an
encoding was proposed by Woods~\cite{woods82} and independently by
Tamassia and Tollis~\cite{tamassia86}. Bipolar orientations of $M$ are
one of the structures which cannot be encoded as $\alpha$-orientations
on $M$, an auxiliary map $M'$ (the angle graph of $M$) has to be used
instead. For Schnyder woods on 3-connected planar maps as well as
bipartite $f$-factors and spanning trees Felsner~\cite{felsner04}
describes encodings as $\alpha$-orientations. He also proves that
  the set of $\alpha$-orientations of a planar map $M$ can always be
  endowed with the structure of a distributive lattice.
This structure on the set of $\alpha$-orientations found applications in
drawing algorithms in~\cite{bfm-cd3cpg-04}, \cite{f-tsotwatsld-05},
and for enumeration and random sampling of graphs
in~\cite{fps-datwatomeatrs-05}.  

Given the existence of a combinatorial structure on a class $\MM_n$ of
planar maps with $n$ vertices, one of the questions of interest is how
many such structures there are for a given map
$M\in\MM_n$. Especially, one is interested in the minimum and maximum
that this number attains on the maps from $\MM_n$. This question has
been treated quite successfully for spanning trees and bipartite
perfect matchings. For spanning trees the Kirchhoff Matrix Tree
Theorem allows to bound the maximum number of
spanning trees of a planar graph with $n$ vertices between~$5.02^n$
and~$5.34^n$, see~\cite{r-tnosptiapg-05,r-phd-05}. Pfaffian
  orientations can be used 
to efficiently calculate the number of bipartite perfect matchings in
the planar case, see for example~\cite{lovaszplummer}.  Kasteleyn has
shown, that the $k\times\ell$ square grid has asymptotically $e^{0.29\cdot
k\ell}\approx 1.34^{k\ell}$ perfect matchings. The number of Eulerian
orientations is studied in statistical physics under the name of ice
models, see~\cite{b-esmism-82} for an overview. In particular
Lieb~\cite{l-treosi-67} has shown that the square grid on the torus
has asymptotically $(8\sqrt{3}/9)^{k\ell}\approx1.53^{k\ell}$ Eulerian
orientations and Baxter~\cite{b-fmtl-69} has worked out the
asymptotics for the triangular grid on the torus as
$(3\sqrt{3}/2)^{k\ell}\approx2.598^{k\ell}$. 

In many cases it is relatively easy to
see which maps in a class $\MM_n$ carry a unique object of a certain
type, while the question about the maximum number is rather
intricate. Therefore, we focus on finding the asymptotics for the
maximum number of $\alpha$-orientations that a map from $\MM_n$ can
carry. The next table gives an overview of the results of this paper for
different instances of $\MM_n$ and $\alpha$. The entry $c$ in the
``Upper Bound'' column is to be read as $O(c^n)$, in the ``Lower
Bound'' column as $\Omega(c^n)$ and for the ``$\approx c$'' entries the
asymptotics are known.

\begin{table}[h]\label{tab:overview}

\begin{tabular}{|l||c|c|}
\hline
 \pbox{.8}{8.2}{\cent{.1}{\bf{Graph class and orientation type}}} &
 \pbox{.8}{3}{\cent{.1}{\bf{Lower bound}}} &
 \pbox{.8}{3}{\cent{.1}{\bf{Upper bound}}}\\\hline
\pbox{.4}{8}{\vspace{.15cm}$\alpha$-orientations on planar maps} &
\pbox{.4}{3}{\cent{.1}{2.598}} &\pbox{.4}{3}{\cent{.1}{3.73}}\\\hline 

\pbox{.4}{8}{\vspace{.15cm}Eulerian orientations} &
\pbox{.4}{3}{\cent{.1}{2.598}} &\pbox{.4}{3}{\cent{.1}{3.73}}\\\hline 

\pbox{.4}{8}{\vspace{.15cm}Schnyder woods on triangulations} &
\pbox{.4}{3}{\cent{.1}{2.37}} &\pbox{.4}{3}{\cent{.1}{3.56}}\\\hline 
\pbox{.4}{8}{\vspace{.15cm}Schnyder woods on the square grid} &
  \multicolumn{2}{c|}{\pbox{.4}{3}{\cent{.1}{$\approx 3.209$}}} \\\hline 
\pbox{.4}{8}{\vspace{.15cm}Schnyder woods on 3-connected planar maps} &
\pbox{.4}{3}{\cent{.1}{3.209}} &\pbox{.4}{3}{\cent{.1}{8}}\\\hline 

\pbox{.4}{8}{\vspace{.15cm}2-orientations on quadrangulations} &
\pbox{.4}{3}{\cent{.1}{1.53}} &\pbox{.4}{3}{\cent{.1}{1.91}}\\\hline 
\pbox{.4}{8}{\vspace{.15cm}bipolar orientations on stacked triangulations} &
  \multicolumn{2}{c|}{\pbox{.4}{3}{\cent{.1}{$\approx 2$}}}
  \\\hline 
\pbox{.4}{8}{\vspace{.15cm}bipolar orientations on outerplanar maps} &
  \multicolumn{2}{c|}{\pbox{.4}{3}{\cent{.1}{$\approx 1.618$}}}
  \\\hline 

\pbox{.4}{8}{\vspace{.15cm}bipolar orientations on the square grid} &
\pbox{.4}{3}{\cent{.1}{2.18}} &\pbox{.4}{3}{\cent{.1}{2.62}}\\\hline 

\pbox{.4}{8}{\vspace{.15cm}bipolar orientations on planar maps} &
\pbox{.4}{3}{\cent{.1}{2.91}} &\pbox{.4}{3}{\cent{.1}{3.97}}\\\hline
\end{tabular}
\end{table}

The paper is organized as follows. In Section~\ref{sec:alpha} we treat
the most general case, where $\MM_n$ is the class of all planar maps
with $n$ vertices and $\alpha$ can be any integer valued
function. We prove an upper bound which applies for every
map and every $\aa$.  In Section~\ref{sec:eulerian} we deal with
Eulerian orientations.  In 
Section~\ref{sec:swont} we consider Schnyder woods on plane
triangulations and in Section~\ref{sec:swong} the more general case of
Schnyder woods on 3-connected planar maps.  We split the treatment of
Schnyder woods because the more direct encoding of Schnyder woods on
triangulations as $\alpha$-orientations yields stronger bounds. In
Section~\ref{sec:swong} we also discuss the asymptotic number of
Schnyder woods on the square grid. Section~\ref{sec:2-or} is dedicated
to 2-orientations of quadrangulations. In
Section~\ref{sec:bipo}, we study bipolar orientations on the square
grid, stacked triangulations, outerplanar maps and planar maps. The upper
bound for planar maps relies on a new encoding of bipolar orientations
of inner triangulations.  
In Section~\ref{subsec:complexity} we discuss the complexity of
counting $\aa$-orientations. In Section~\ref{subsec:perfmatch} we show
how counting $\alpha$-orientations can be reduced to 
counting (not necessarily planar) bipartite perfect matchings and the
consequences of this connection are discussed as well. We conclude
with some open problems.


\section{Counting $\alpha$-Orientations}\label{sec:alpha}


 A \term{planar map} $M$ is a simple planar graph 
 $G$ together with a fixed crossing-free embedding of $G$ in the Euclidean
 plane. In particular, $M$ has a designated outer (unbounded) face. We
 denote the sets of vertices, edges and faces of a given planar map by
 $V$, $E$, and $\mathcal{F}$, and their respective cardinalities  by
 $n$, $m$ and $f$. The degree of a vertex $v$ will be denoted by 
 $d(v)$.
 
Let $M$ be a planar map and $\alpha: V \rightarrow \NN$. An orientation
 $X$ of the edges of $M$ is an \term{$\alpha$-orientation} if for all $v\in V$ exactly
 $\aa(v)$ edges are directed away from $v$ in $X$.

 Let $X$ be an $\alpha$-orientation of $M$ and let $C$ be a directed cycle
 in $X$. Define $X^C$ as the orientation obtained from $X$ by reversing all
 edges of $C$. Since the reversal of a directed cycle does not affect
 out-degrees the orientation $X^C$ is also an
 $\alpha$-orientation of $M$. The plane embedding of $M$ allows us to
 classify a directed simple cycle as clockwise (\term{cw-cycle}) if the
 interior, $\Int(C)$, is to the right of $C$ or as counterclockwise
 (\term{ccw-cycle}) if $\Int(C)$ is to the left of $C$. If $C$ is a
 ccw-cycle of $X$ then we say that $X^C$ is \term{left of} $X$ and $X$
 is \term{right of} $X^C$. Felsner proved the following theorem in~\cite{felsner04}.  
\begin{theorem}\label{thm:d-lat}
   Let $M$ be a planar map and $\alpha: V \to \NN$. The set of
   $\alpha$-orientations of $M$ endowed with the transitive closure of the
   `left of' relation is a distributive lattice. 
 \end{theorem}
The following observation is easy, but useful. Let $M$ and 
$\aa:V\rightarrow \NN$ be given,  $W\subset V$ and $E_W$ the edges of
$M$ with one endpoint in $W$ and the other endpoint in $V\setminus
W$. Suppose all edges of $E_W$ are directed away from $W$ in some
$\alpha$-orientation $X_0$ of $M$. The demand of $W$ for $\sum_{w\in
  W}\aa(w)$ outgoing edges forces all edges 
in $E_W$ to be directed away from $W$ in every $\alpha$-orientation
of~$M$. Such an edge with the same direction in every $\alpha$-orientation
is a \term{rigid edge}.

We denote the number of $\alpha$-orientations of $M$ by
$r_{\alpha}(M)$. Let $\MM$ be a family of  pairs $(M,\aa)$ of a planar
map and an out-degree function.
Most of this paper is concerned with lower and upper
bounds for $\max_{(M,\aa)\in\mathcal{M}}r_{\alpha_M}(M)$ for some family
$\MM$. In Section~\ref{sec:genupperbound}, we deal with
bounds which apply to all $M$ and $\alpha$, while later sections will
be concerned with special instances.


\subsection{An Upper Bound for the Number $\aa$-Orientations}\label{sec:genupperbound}


A trivial upper bound for the number of $\aa$-orientations on $M$ is
$2^m$ as any edge can be directed in two ways. The following easy but
useful lemma improves the trivial bound. 

\begin{lemma}\label{lem:tree} Let $M$ be a planar map, $A\subset E$
  a cycle free subset of edges of $M$, and $\aa$ a function $\alpha:V\rightarrow \NN$. Then, there are
  at most $2^{m-|A|}$ $\alpha$-orientations of $M$. Furthermore, $M$ has less than $4^n$
  $\aa$-orientations.  
\end{lemma}
\proof 
Let $X$ be an arbitrary but fixed orientation out of the
$2^{m-|A|}$ orientations of the edges of 
$E\setminus A$. It suffices to show that $X$ can be extended to an
$\alpha$-orientation of $M$ in at most one way. 
We proceed by induction on $|A|$. The base case $|A|=0$ is
trivial. If $|A|>0$, then, as $A$ is cycle free, there is a vertex $v$,
which is incident to exactly one edge $e\in A$. If $v$ has out-degree
$\alpha(v)$ respectively $\alpha(v)-1$ in $X$, then $e$ must be
directed towards 
respectively away from $v$. In either case the direction of $e$ is
determined by $X$, and by induction there is at most one way to
extend the resulting orientation of $E\setminus(A-e)$ to an
$\aa$-orientation of $M$. If $v$ does not have out-degree $\aa(v)$ or
$\aa(v)-1$ in $X$, then there is no extension of $X$ to an
$\aa$-orientation of $M$. The bound $2^{m-n+1}<4^n$ follows by
choosing $A$ to be a spanning forest and applying  Euler's formula. \qed

A better upper bound for general $M$ and $\alpha$ will be given in
Proposition~\ref{prop:generalupperbound}. The following lemma is
needed for the proof.
\begin{lemma}\label{lem:edgeinsert} 
 Let $M$ be a planar map with $n$ vertices that has an
  independent set $I_2$ of $n_2$ vertices which have degree 2 in $M$. Then,
  $M$ has at most $(3n-6) -(n_2 - 1)$ edges.
\end{lemma}
\proof 
Consider a triangulation $T$ extending $M$ and let $B$ be the set of
additional edges, i.e., of edges of $T$ which are not in $M$.  If $n=3$ the
conclusion of the lemma is true and we may thus assume $n
>3$ for the rest of the proof. Hence, there are no vertices of degree 2
in $T$, and every vertex of $I_2$ must be incident to at least one edge
from $B$. If there is a vertex $v\in 
I_2$, which is incident to exactly one edge from $B$, then $v$ and its
incident edges can be deleted from $I_2$, from $M$ and from $T$,
whereby the result follows by induction. The last case is that all
vertices of $I_2$ have at least two incident edges in 
$B$.  Since every edge in $B$ is incident to at most two vertices from
$I_2$ it follows that $|I_2|\leq |B|$. Therefore, $|E_M| = |E_T| - |B|
\leq |E_T| - |I_2| = (3n-6) - n_2$.
\qed

\ni\textbf{Remark.} It can be seen from the above proof, that
$K_{2,n_2}$ plus the edge between the two vertices of degree $n_2$ is
the unique graph to which only $n_2-1$ edges can be added. For
every other graph at least $n_2$ edges can be added.\medskip

\begin{proposition}\label{prop:generalupperbound}
Let $M$ be a planar map, $\alpha:V\rightarrow\NN$, and $I=I_1\cup I_2$ an
independent set of $M$, where $I_2$ is the subset of degree 2 vertices
in $I$. Then, $M$ has at most 
\begin{equation}\label{eq:generalupperbound}2^{2n-4-|I_2|}\cdot\prod_{v\in I_1}
\left(\frac{1}{2^{d(v)-1}}\bin{d(v)}{\aa(v)}\right)\end{equation}
$\alpha$-orientations.
\end{proposition}

\proof We may assume that $M$ is connected. Let
$M_i$, for $i=1,\ldots c$, be the components of $M-I$. We claim that
$M$ has at most $(3n-6) - (c-1) - (|I_2|-1)$ edges. Note, that every
component $C$ of $M-I$ must be connected to some other component $C'$
via a vertex $v\in I$ such that the edges $vw$ and $vw'$ with $w\in C$
and $w'\in C'$ form an angle at $v$. As $w$ and $w'$ are in different
connected components the edge $ww'$ is not in $M$ and we can add it
without destroying planarity.  We can add at least $c-1$ edges not
incident to $I$ in this fashion. Thus, by Lemma~\ref{lem:edgeinsert}
we have that $m+(c-1)\leq 3n-6-(I_2-1)$.

Let $S'$ be a spanning forest of $M-I$, and
let $S$ be obtained from $S'$ by adding one edge incident to every
$v\in I$. Then, $S$ is a forest with $n-c$ edges. By Lemma~\ref{lem:tree} $M$ has at most
$2^{m-|S|}$ $\aa$-orientations and by Lemma~\ref{lem:edgeinsert} 
\[m-|S|\leq (3n-6)-(c-1)-(|I_2|-1)-(n-c)=2n-4-|I_2|.\]

For every vertex $v\in I_1$ there are $2^{d(v)-1}$ possible
orientations of the edges of 
$M-S$ at $v$. Only the orientations with $\aa(v)$ or $\aa(v)-1$
outgoing edges at $v$ can potentially be completed to an
$\alpha$-orientation of $M$. Since $I_1$ is an independent set it
follows that $M$ has at most 
\begin{equation}\label{eqn:generalbound}
  2^{m-|S|}\cdot\prod_{v\in I_1} 
  \frac{1}{2^{d(v)-1}}\left(\bin{d(v)-1}{\aa(v)}+\bin{d(v)-1}{\aa(v)-1}\right)\leq 
2^{2n-4-|I_2|}\cdot\prod_{v\in I_1}
\left(\frac{1}{2^{d(v)-1}}\bin{d(v)}{\aa(v)}\right)
\end{equation}
$\aa$-orientations. \qed

\begin{corollary}\label{cor:alphauppperbound}
Let $M$ be a planar map and $\alpha:V\rightarrow\NN$. Then, $M$ has at
most $3.73^n$ $\alpha$-orientations.
\end{corollary}
\proof Since $M$ is planar the Four Color Theorem implies, that it has an 
independent set $I$ of size $|I|\geq n/4$. Let $I_1$, $I_2$ be as
above. Note, that for $d(v)\geq 3$ 
\begin{equation}\label{eqn:upperbound} \frac{1}{2^{d(v)-1}}\bin{d(v)}{\aa(v)}\leq
\frac{1}{2^{d(v)-1}}\bin{d(v)}{\lfloor d(v)/2\rfloor}\leq \frac{3}{4}.
\end{equation}
Thus, the result follows from
Proposition~\ref{prop:generalupperbound}, as 
$$
2^{2n-4-|I_2|}\left(\frac{3}{4}\right)^{|I_1|}\leq
2^{2n-4}\left(\frac{3}{4}\right)^\frac{n}{4} 
\leq 3.73^n.
$$ \qed

\ni\textbf{Remark.} The best lower bound for general $\alpha$ and $M$,
which we can prove, comes from Eulerian orientations of the triangular grid, see
Section~\ref{sec:eulerian}. 


\subsection{Grid Graphs}\label{sec:grids}


Enumeration and counting of different combinatorial
structures on grid graphs have received a lot of attention in the
literature, see~e.g.~\cite{b-esmism-82,calkinwilf,l-treosi-67}. 
In Section~\ref{sec:eulerian} we present a family of graphs that have
asymptotically at least $2.598^n$ Eulerian orientations. This family
is closely related to the grid graph, and throughout the paper we will
use different relatives of the grid graph to obtain lower bounds. We collect the
definitions of these related families here. 

The grid graph $G_{k,\ell}$ with $k$ rows and $\ell$ columns is defined
as follows. The vertex set is 
\[V_{k,\ell}=\{(i,j)\mid 1\leq i\leq k,\;1\leq j\leq \ell\}.\]
The edge set $E_{k,\ell}=E^H_{k,\ell}\cup E^V_{k,\ell}$ consists of
horizontal edges 
\[
E^H_{k,\ell}=\Big\{\{(i,j),(i,j+1)\}\mid 1\leq i\leq k,\;1\leq j\leq\ell-1\Big\}
\]
and vertical edges
\[
E^V_{k,\ell}=\Big\{\{(i,j),(i+1,j)\}\mid 1\leq i\leq k-1,\; 1\leq j\leq\ell\Big\}.
\]
We denote the $i$th vertex row by $V^R_i=\{(i,j)\mid 1\leq j\leq
\ell\}$ and the $j$th vertex column by $V^C_j=\{(i,j)\mid
1\leq i\leq k\}$. The $j$th edge column $E^C_j$ is defined as
$E^C_j=\{\{(i,j),(i,j+1)\}\mid 1\leq i\leq k\}$. 
The number of bipolar orientations of $G_{k,\ell}$ is studied in
Section~\ref{sec:grid}.

The grid on the torus $G^T_{k,\ell}$ is obtained from
$G_{k+1,\ell+1}$ by identifying $(1,i)$ and $(k+1,i)$ as well as
$(j,1)$ and $(j,\ell+1)$ for all $i$ and $j$, see
Figures~\ref{fig:squaregrids3}~(a) and~(b).  Edges of the form $\{(i,1),(i,\ell)\}$
are called horizontal wrap-around edges while those of the form
$\{(1,j),(k,j)\}$ are the vertical wrap-around edges. Note that
$G_{k,\ell}$ can be obtained from $G^T_{k,\ell}$ by deleting the $k$
horizontal and the $\ell$ vertical wrap-around edges. 

Lieb~\cite{l-treosi-67} shows that $G^T_{k,\ell}$ has asymptotically
$(8\sqrt{3}/9)^{k\ell}$ Eulerian orientations. His analysis involves
the calculation of the dominant eigenvalue of a so-called transfer
matrix, see also Section~\ref{sec:2-or}.

{\def\SetFigFont#1#2#3#4#5{%
  \fontsize{#1}{#2pt}%
  \fontfamily{#3}\fontseries{#4}\fontshape{#5}%
  \selectfont}%

\PsFigCap{47}{squaregrids3}
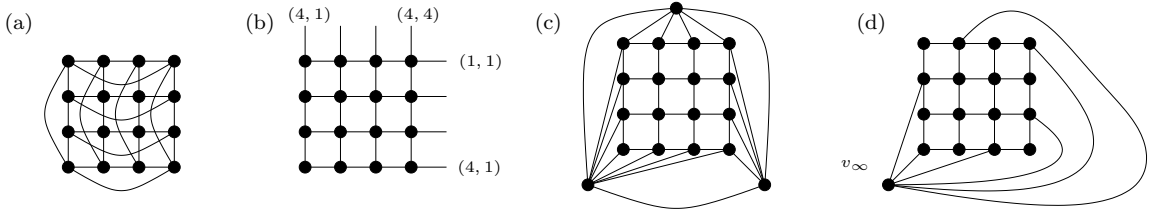
{Two illustrations of $G^T_{4,4}$, the augmented grid $G^*_{4,4}$, and the quadrangulation $G^\Box_{4,4}$.}
}

We consider the number of Schnyder woods on the
augmented grid $G^*_{k,\ell}$ in Section~\ref{sec:swong}, see
Figure~\ref{fig:squaregrids3}~(c). The augmented
grid is obtained from $G_{k,\ell}$ by adding a triangle with vertices
$\{a_1,a_2,a_3\}$ to the outer face. The triangle is connected to the
boundary vertices of the grid as follows. The vertex $a_1$ is adjacent
to all vertices of $V^R_1$, $a_2$ is adjacent the vertices from
$V^C_\ell$ and $a_3$ to the vertices from $V^R_k\cup V^C_1$. 

When we consider 2-orientations in Section~\ref{sec:2-or} we use the
quadrangulation $G^\Box_{k,\ell}$, see
Figure~\ref{fig:squaregrids3}~(d). It is obtained from the grid   
$G_{k,\ell}$ by adding one vertex $v_\infty$ to the outer face which
is adjacent to every other vertex of the boundary such that $(1,1)$ is
not adjacent to $v_\infty$. For $k$
and $\ell$ even this graph is closely related to the torus grid
$G_{k,\ell}^T$, which can be obtained from $G^\Box_{k,\ell}$ by
reassigning end vertices of edge as follows.
\[
\begin{array}{llcll} 
  \{(1,j),v_\infty\}\rightarrow\{(1,j),(k,j)\} & 2\leq j\leq \ell & \hspace{.5cm}
& \{(k,j),v_\infty\}\rightarrow\{(k,j),(1,j)\} & 2\leq j\leq \ell\\
  \{(i,1),v_\infty\}\rightarrow\{(i,1),(i,\ell)\} & 2\leq i\leq k &
& \{(i,\ell),v_\infty\}\rightarrow\{(i,\ell),(i,1)\} & 2\leq i\leq k \\
\end{array}
\]

Since $k,\ell$ are even this does not
create parallel edges and the resulting graph is $G^T_{k,\ell}$
minus the edges $e_1=\{(1,1),(1,\ell)\}$ and
$e_2=\{(1,1),(k,1)\}$.

We also use the triangular grid $T_{k,\ell}$ in
Sections~\ref{sec:eulerian} and~\ref{bipolar}. It
is obtained from
$G_{k,\ell}$ by adding the diagonal edges $\{(i,j),(i-1,j+1)\}$ for 
$2\leq i\leq k$ and $1\leq j\leq \ell-1$, see Figure~\ref{fig:trigrids}~(a). The
augmented triangular grid $T^*_{k,\ell}$, which we need in
Section~\ref{sec:swont} is obtained in the same way from
$G^*_{k,\ell}$, see Figure~\ref{fig:trigrid2}. 

The terms vertex row, vertex column and edge column are used for the
triangular grid analogously to the definition above for $G_{k,\ell}$. 

{\def\SetFigFont#1#2#3#4#5{%
  \fontsize{#1}{#2pt}%
  \fontfamily{#3}\fontseries{#4}\fontshape{#5}%
  \selectfont}%

\PsFigCap{30}{trigrids}
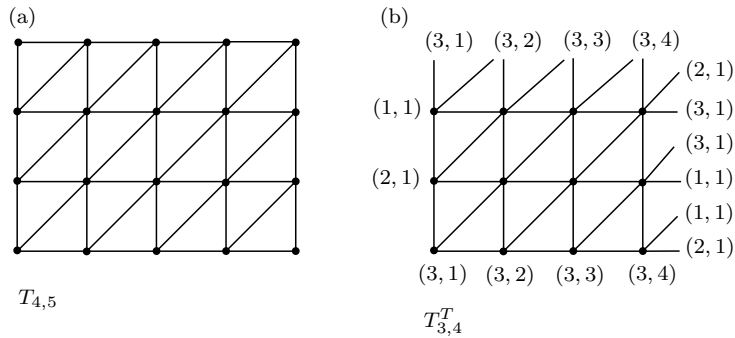
{The triangular grid $T_{4,5}$, and $T^T_{3,4}$.}
}

We also use the triangular grid on the torus $T^T_{k,\ell}$, see
Figure~\ref{fig:trigrids}~(b). We adopt the 
definition from~\cite{b-fmtl-69}, therefore it differs slightly from
that of the square grid on the torus. More precisely, instead of
identifying vertices $(i,\ell+1)$ and $(i,1)$ we identify vertices
$(i,\ell+1)$ and $(i-1,1)$ (and $(1,\ell+1)$ with $(k,1)$) to obtain 
$T^T_{k,\ell}$ from $T_{k,\ell}$. This boundary condition is called
helical. The wrap-around edges are defined analogously to the square
grid case. 

Baxter~\cite{b-fmtl-69} was able to determine the  
exponential growth factor of Eulerian orientations of $T^T_{k,\ell}$ as
$k,\ell\rightarrow\infty$. Baxter's analysis uses similar techniques
as Lieb's~\cite{l-treosi-67} and yields an asymptotic growth rate of
$(3\sqrt{3}/2)^{k\ell}$.


\subsection{A Lower Bound Using Eulerian Orientations}\label{sec:eulerian}

%

Let $M$ be a planar map such that every $v\in V$ has even
degree and let $\aa$ be defined as $\aa(v)=d(v)/2,\;\forall v\in
V$. The corresponding $\aa$-orientations of $M$ are known as
\term{Eulerian orientations}. Eulerian orientations are
exactly the orientations which maximize the binomial coefficients in 
equation (\ref{eq:generalupperbound}). The lower bound in the next
theorem is the best lower bound we have for $\max_{(M,\aa)\in\MM} r_{\aa_M}(M)$,
where $\MM$ is the set of all planar maps and no restrictions are made
for $\aa$.
\begin{theorem}\label{thm:eulerian}
Let $\mathcal{M}_n$ denote the set of all planar maps with
$n$  vertices and  $\mathcal{E}(M)$ the set of Eulerian orientations of
$M\in\mathcal{M}_n$. Then, for $n$ big enough, \[2.59^n\leq(3\sqrt{3}/2)^{k\ell}\leq
\max_{M\in\mathcal{M}_n}|\mathcal{E}(M)|\leq3.73^{n}.\]
\end{theorem}
\proof The upper bound is the one from Corollary~\ref{cor:alphauppperbound}.
For the lower bound consider the triangular torus grid $T^T_{k,\ell}$.
 As mentioned above Baxter~\cite{b-fmtl-69} was able to determine the  
exponential growth factor of Eulerian orientations of $T^T_{k,\ell}$ as
$k,\ell\rightarrow\infty$. Baxter's analysis uses
eigenvector calculations and yields an asymptotic growth rate of
$(3\sqrt{3}/2)^{k\ell}$. This graph can be made into a planar map
$T_{k,\ell}^+$ by introducing a new vertex $v_\infty$ which is
incident to all the wrap-around edges. This way all crossings between
wrap-around edges can be substituted by $v_\infty$. As every Eulerian
orientation of  $T^T_{k,\ell}$ yields a Eulerian orientation of  
$T_{k,\ell}^+$ this graph has at least $(3\sqrt{3}/2)^{k\ell}\geq
2.598^{k\ell}$ Eulerian orientations for $k,\ell$ big enough. \qed


\section{Counting Schnyder Woods}\label{sec:sw}

Schnyder woods for triangulations have been introduced as a tool for
graph drawing and graph dimension theory
in~\cite{schnyder89,schnyder90}. Schnyder woods for 3-connected 
planar maps are introduced in~\cite{felsner01}. Here we review the
definition of Schnyder woods and explain how they are encoded as
$\alpha$-orientations. For a comprehensive introduction see
e.g.~\cite{f-gga-04}. 

Let $M$ be a planar map with three vertices $a_1,a_2,a_3$  occurring 
in clockwise order on the outer face of $M$. A suspension $\MS$ of $M$ is
obtained by attaching a half-edge that reaches into the outer face to
each of these \term{special vertices}.

Let $\MS$ be a suspended 3-connected planar map.
A \term{Schnyder wood} rooted at $a_1,a_2,a_3$ is an orientation and coloring
of the edges of $\MS$ with the colors $1,2,3$
satisfying the following rules.

\Item(W1)
Every edge $e$ is oriented in one direction or in two opposite directions.
The directions of edges are colored such that
if $e$ is bidirected the two directions have
distinct colors.

\Item(W2)
The half-edge at $a_i$ is directed outwards and has color $i$.

\Item(W3)
Every vertex $v$ has out-degree one in each color.
The edges $e_1,e_2,e_3$ leaving $v$ in colors $1,2,3$ occur in clockwise
order. Each edge entering $v$ in color $i$ enters $v$ in the clockwise
sector from $e_{i+1}$ to $e_{i-1}$, see Figure~\ref{fig:orient+notunique}~(a).

\Item(W4)
There is no interior face the boundary of which is a monochromatic directed cycle.
\medskip

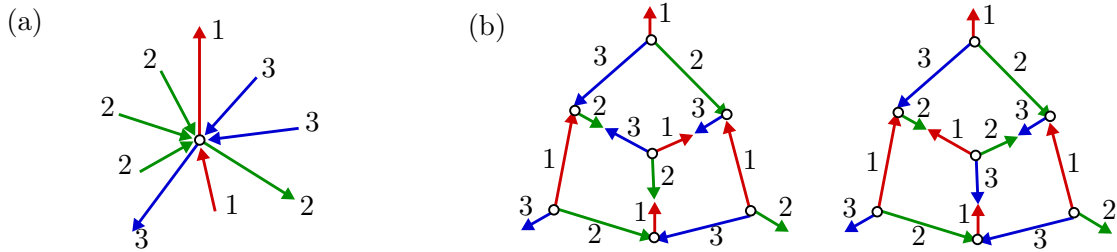
\PsFigCap{80}{orient+notunique}{The left part shows edge orientations
  and edge colors at a vertex, the right part two different Schnyder
  woods with the same underlying orientation.}

In the context of this paper the choice of the suspension vertices is
not important and we refer to the Schnyder wood of a planar map,
without specifying the suspension explicitly. 

Let $\MS$ be a planar map with a Schnyder wood. Let $T_i$ denote the
digraph induced by the directed edges of color $i$.  Every inner
vertex has out-degree one in $T_i$ and in fact $T_i$ is a directed spanning
tree of $M$ with root $a_i$. 

In a Schnyder wood on a triangulation only the three
outer edges are bidirected. This is because the three spanning trees 
have to cover all $3n-6$ edges of the 
triangulation and the edges of the outer triangle must be bidirected
because of the rule of vertices. Theorem~\ref{thm:swont<->or} says,
that the edge orientations together with the colors of the special vertices
are sufficient to encode a Schnyder wood on a triangulation,
the edge colors can be deduced, for a proof see~\cite{fraysseixmendez01}. 

\begin{theorem}\label{thm:swont<->or}
Let $T$ be a plane triangulation, with vertices $a_1,a_2,a_3$ occuring
in clockwise order on the outer face. Let
$\alpha_T(v):=3$ if $v$ is an internal vertex and $\alpha_T(a_i):=0$
for $i=1,2,3$. Then, there is a bijection 
between the Schnyder woods of $T$ and the $\alpha_T$-orientations of
the inner edges of $T$.
\end{theorem}

In the sequel we refer to an $\alpha_T$-orientation
simply as a 3-orientation. Schnyder woods on 3-connected planar maps
are in general not uniquely determined by the edge
orientations, see Figure~\ref{fig:orient+notunique}~(b). Nevertheless, there is
a bijection between the Schnyder woods
of a 3-connected planar map $M$ and certain $\aa$-orientations on
a related planar map $\MT$, see~\cite{felsner04}. 

In order to describe the bijection precisely, we first define the
\term{suspension dual} $\MSD$ of 
$\MS$, which is obtained from the dual $M^*$ of $M$ as follows. Replace
the vertex $v^*_\infty$, which represents the unbounded face of $M$ in
$M^*$, by a triangle on three new vertices
$b_1,b_2,b_3$.  Let $P_i$ be the path from  $a_{i-1}$ to $a_{i+1}$
on the outer face of $M$ which avoids $a_i$. In $\MSD$ the edges dual
to those on $P_i$ are incident to $b_i$  instead of $v^*_\infty$. Adding a ray
 to each of the $b_i$ yields $\MSD$. An example is  
given in Figure~\ref{fig:PrimalDualSW3}.

\PsFigCap{22}{PrimalDualSW3}
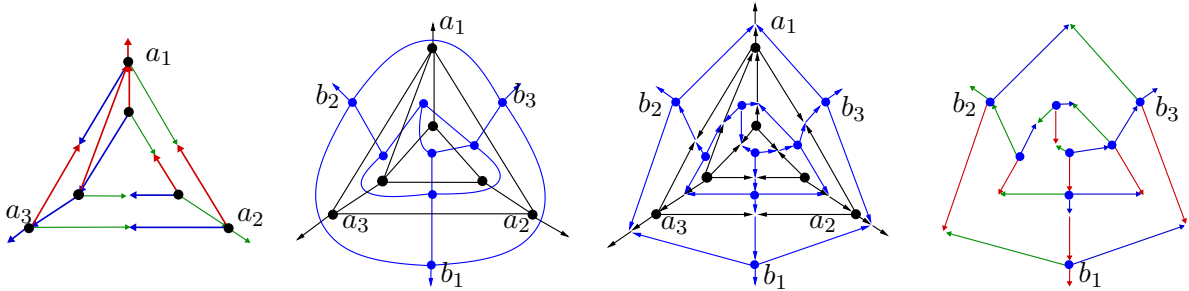
{A Schnyder wood, the primal and the dual graph, the oriented primal
  dual completion and the dual Schnyder wood.}
\begin{proposition}\label{prop:swp<->swd}
Let $\MS$ be a suspended planar map.
There is a bijection between the Schnyder woods of $\MS$ and 
the Schnyder woods of the suspension dual
$\MSD$. Figure~\ref{fig:primaldualedges} illustrates how the coloring
and orientation of a pair of a primal and a dual edge are related.  
\end{proposition}
%
\PsFigCap{30}{primaldualedges}
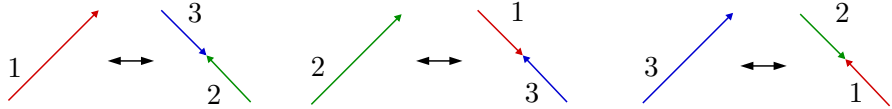
{The three possible oriented colorings of a pair of a primal and a dual edge.}

%
%
The completion $\MT$ of $\MS$ and $\MSD$ is obtained by superimposing
the two graphs such that exactly the primal dual pairs of edges
cross, see Figure~\ref{fig:PrimalDualSW3}. In the completion $\MT$ the
common subdivision of each crossing 
pair of edges is replaced by a new edge-vertex. Note
that the rays emanating from the three special vertices of $\MS$  cross
the three edges of the triangle induced by $b_1,b_2,b_3$ and thus produce edge
vertices. The six rays emanating into the unbounded face of the
completion end at a new vertex $v_\infty$ placed in this unbounded
face. A pair of corresponding Schnyder woods on $\MS$ and $\MSD$
induces an orientation of $\MT$ which is an $\alpha_S$-orientation where
\[
\alpha_S(v) =\left\{\begin{array}{ll} 3&\mbox{ for primal and dual vertices}\\
  1& \mbox{ for edge vertices}\\
  0& \mbox{ for }v_\infty.\end{array}\right.
\]
Note, that a pair of a primal and a dual edge always consists of a
unidirected and a bidirected edge, which explains why $\alpha_S(v_e) =
1$ is the right choice. Theorem~\ref{thm:swong<->or} says,
that the edge orientations of $\MT$ are sufficient to encode a Schnyder
wood of $\MS$, the edge colors can be deduced, for a proof see
\cite{felsner04}. 

\begin{theorem}\label{thm:swong<->or}
The Schnyder woods of a suspended planar map $\MS$ are in bijection with the
$\alpha_S$-orientations of $\MT$.
\end{theorem}
 
In the rest of this section we give asymptotic bounds for the maximum number of
Schnyder woods on planar triangulations and 3-connected planar
maps. We treat these two 
classes separately because the more direct bijection from
Theorem~\ref{thm:swont<->or} allows us to obtain a better upper bound
for Schnyder woods on triangulations than for the general case. We also
have a better lower bound for the general case of Schnyder woods on
3-connected planar maps than for the restriction to
triangulations. 

\term{Stacked triangulations} are plane triangulations which can be obtained
from a triangle by iteratively adding vertices of degree 3 into
bounded faces. The stacked triangulations are
exactly the plane triangulations which have a unique Schnyder wood and
we have a generalization of this well-known result for general
3-connected planar maps, which we state here without a proof.
\begin{theorem}\label{thm:unique}
All 3-connected planar maps, which have a unique
Schnyder wood, can be constructed from the
unique Schnyder wood on the triangle by the six operations show in
Figure~\ref{fig:reductions} read from left to right. 
\end{theorem}
\PsFigCap{25}{reductions}
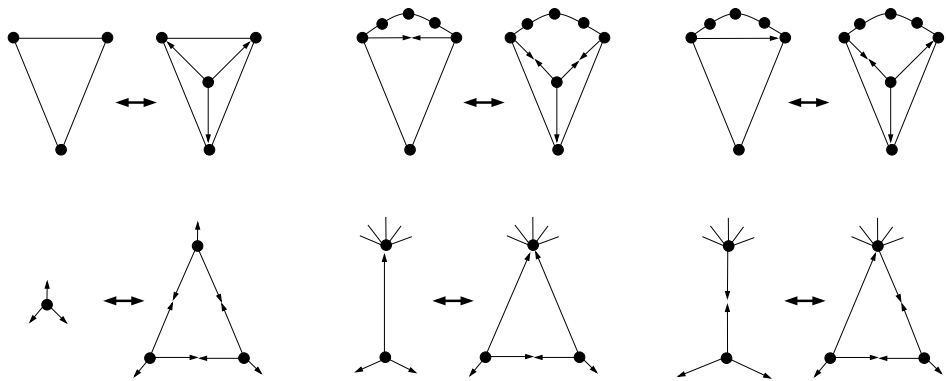
{Using the three primal operations in the first row and their duals in the
  second row every graph with a unique Schnyder wood can be constructed.} 

\subsection{Schnyder Woods on Triangulations}\label{sec:swont}

Bonichon~\cite{b-brmpgpdp-02} found a bijection between Schnyder woods
on triangulations with $n$ vertices and pairs of non-crossing
Dyck-paths, which implies that there are $C_{n+2}C_n-C_{n+1}^2$
Schnyder woods on triangulations with $n$ vertices. By $C_n$
we denote the $n$th Catalan number $C_n=\bin{2n}{n}/(n+1)$. Hence,
asymptotically there are about 
$16^n$ Schnyder woods on triangulations with $n$ vertices. Tutte's
classic result~\cite{t-acopt-62} yields that there are 
asymptotically about $9.48^n$ plane triangulations on $n$
vertices. See~\cite{ps-ocasot-04} for a proof of Tutte's formula
using Schnyder woods. The two results together imply that a
triangulation with $n$ vertices has on average about $1.68^n$ Schnyder woods. 
The next theorem is concerned with the maximum number of
Schnyder woods on a fixed triangulation.

\begin{theorem}\label{thm:swontbounds} Let $\mathcal{T}_n$ denote the
  set of all plane triangulations with $n$  vertices and
  $\mathcal{S}(T)$ the set of Schnyder woods of $T\in\mathcal{T}_n$. Then, 
  \[2.37^n\leq \max_{T\in\mathcal{T}_n}|\mathcal{S}(T)|\leq3.56^{n}.\]
\end{theorem}
The upper bound follows from
Proposition~\ref{prop:generalupperbound} by using that for $d(v)\geq
3$ $$\bin{d(v)}{3}\cdot2^{1-d(v)}\leq \frac{5}{8}.$$ 
For the proof of the lower bound we use the augmented triangular grid
$T^*_{k,\ell}$. Figure~\ref{fig:trigrid2} shows a canonical Schnyder 
wood on $T^*_{k,\ell}$ in which the vertical edges are directed up, the
horizontal edges to the right and diagonal ones left-down.

%
{\def\SetFigFont#1#2#3#4#5{%
  \fontsize{#1}{#2pt}%
  \fontfamily{#3}\fontseries{#4}\fontshape{#5}%
  \selectfont}%

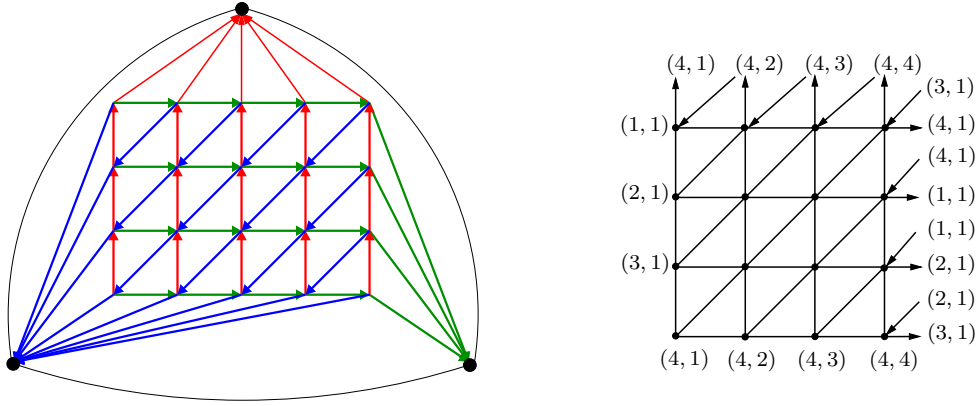
\PsFigCap{30}{trigrid2}{The graphs $T^*_{4,5}$ with a canonical
  Schnyder wood and  $T_{4,4}$ with the additional edges simulating
  Baxter's boundary conditions.}
}
%
Instead of working with the 3-orientations of $T^*_{k,\ell}$ we use
$\alpha^*$-orientations of $T_{k,\ell}$ where
\[\alpha^*(i,j)=\left\{
  \begin{array}{ll}
    3 & \mbox{ if $2\leq i\leq k-1$ and $2\leq i\leq \ell-1$}\\
    1 & \mbox{ if } (i,j)\in\{(1,1),(1,\ell),(k,\ell)\}\\
    2 & \mbox{ otherwise.}
  \end{array}\right.\]
For the sake of simplicity, we refer to $\alpha^*$-orientations of
$T_{k,\ell}$ as $3$-orientations.

Intuitively, $T_{k,\ell}$ promises to be a good
candidate for a lower bound because the canonical orientation shown in
Figure~\ref{fig:trigrid2} on the left has many directed cycles. We
formalize this intuition in the next proposition, which we restrict to
the case $k=\ell$ only to make the notation easier.

\begin{proposition}\label{prop:easylow} The graph $T^*_{k,k}$ has at least
  $2^{5/4(k-1)^2}$ Schnyder woods.  For $k$ big enough $T^*_{k,k}$ has
  \[2.37^{k^2+3}\leq |\mathcal{S}(T^*_{k,k})|\leq2.599^{k^2+3}\]
\end{proposition}
\proof The face boundaries of the triangles of $T_{k,k}$ can be
partitioned into two classes $\mathcal{C}$ and $\mathcal{C}'$ of
directed cycles, such that each class has cardinality $(k-1)^2$ and no two cycles from the
same class share an edge. Thus, a cycle $C\in\mathcal{C}'$ shares an edge with
three cycles from $\mathcal{C}$ if it does not share an edge with the
outer face of $T_{k,k}$ and otherwise it shares an edge with one or two
cycles from  $\mathcal{C}$.

For any subset $D$ of $\mathcal{C}$ reversing all the cycles in $D$
yields a 3-orientation of  $T_{k,k}$, and we can encode this
orientation as a 0-1-sequence of length $(k-1)^2$. After performing
the flips of a given 0-1-sequence $a$, an inner cycle 
$C'\in\mathcal{C}'$ is directed if and only if either all or
none of the three cycles sharing an edge with $C'$ have been
reversed. If $C'\in\mathcal{C}'$ is a boundary cycle, then it is
directed if and only none of the adjacent cycles from $\mathcal{C}$ has been
reversed. Thus the number of different cycle flip sequences on
$\mathcal{C}\cup \mathcal{C'}$ is bounded from below by
\[\sum_{a\in\{0,1\}^{(k-1)^2}}2^{\sum_{C'\in\mathcal{C}'}X_{C'}(a)}.\]
Here $X_{C'}(a)$ is an indicator function which takes value 1 if $C'$
is directed after performing the flips of $a$ and 0 otherwise.

We now assume that every $a\in\{0,1\}^{(k-1)^2}$ is chosen uniformly
at random. The expected value of the above function is then 
\[\EE[2^{\sum X_{C'}}]=
\frac{1}{2^{(k-1)^{2}}}\sum_{a\in\{0,1\}^{(k-1)^2}}
  2^{\sum_{C'\in\mathcal{C}'}X_{C'}(a)}.\]  

Jensen's inequality  $\EE[\varphi(X)]\geq
\varphi(\EE[X])$ holds for a random variable $X$ and a convex
function $\varphi$. Using this we derive that
\[\EE[2^{\sum X_{C'}}]\geq 
2^{\EE[\sum X_{C'}]}=2^{\sum \PP[C'\; \mbox{{\scriptsize flippable}}]}.\]
The probability that $C'$ is flippable is at least $1/4$. For $C'$ which does
not include a boundary edge the probability depends only on
the three cycles from $\mathcal{C}$ that share an edge with $C'$ and
two out of the eight flip vectors for these three cycles make $C'$
flippable. A similar reasoning applies for $C'$ including a boundary
edge. Altogether this yields that 
\[\sum_{a\in\{0,1\}^{(k-1)^2}}2^{\sum_{C'\in\mathcal{C}'}X_{C'}(a)}
\geq 2^{(k-1)^2} \cdot 2^{(k-1)^2/4}.\]
Different cycle flip sequences yield different Schnyder woods. The
orientation of an edge is easily determined. The edge direction  is
reversed with 
respect to the canonical orientation if and only if exactly one of the two
cycles on which it lies has been flipped. We can
tell a flip sequence apart from its complement by looking at the
boundary edges. 

For the upper bound we use Baxter's result for Eulerian orientations 
on the torus $T^T_{k,\ell}$ (see Sections~\ref{sec:grids} and~\ref{sec:eulerian}). Every
3-orientation of $T_{k,\ell}$ plus the wrap-around edges, oriented as shown in 
Figure~\ref{fig:trigrid2} on the right, yields a Eulerian orientation
of $T^T_{k,\ell}$. We deduce that $T^*_{k,\ell}$ has at most  $2.599^n$
Schnyder woods. \qed

\ni\textbf{Remark.} Let us briefly come back to the number of Eulerian
orientations of $T^T_{k,\ell}$, which
was mentioned in Sections~\ref{sec:grids} and~\ref{sec:eulerian} and
in the above proof. 
There are only $2^{2(k+\ell)-1}$ different orientations of the
wrap-around edges. By the pigeon hole principle there is an 
orientation of these edges which can be extended to
a Eulerian orientation of $T^T_{k,\ell}$ in asymptotically
$(3\sqrt{3}/2)^{k\ell}$ ways. Thus, there are out-degree functions
$\alpha_{k\ell}$ for $T_{k,\ell}$ such that there are asymptotically
$2.598^{k\ell}$ 
$\aa_{k\ell}$-orientations. Note, however, that directing all the 
wrap-around edges away from the vertex to which they are attached in
Figure~\ref{fig:trigrid2} induces a unique Eulerian orientation of
$T_{k,\ell}$. 

We have not been able to specify orientations of the
wrap-around edges, which allow to conclude that $T_{k,\ell}$ has
$(3\sqrt{3}/2)^{k\ell}$ 3-orientations with these boundary conditions. In
particular we have no proof that Baxter's result also gives a lower
bound for the number of 3-orientations.  

\subsection{Schnyder Woods on the Grid and 3-Connected Planar
  Maps}\label{sec:swong} 

In this section we discuss bounds on the number of Schnyder
woods on 3-connected planar maps. The lower bound comes from the
grid. The upper bound for this case is much larger than the one for
triangulations. This is due to the encoding of Schnyder woods by
3-orientations  on the primal dual completion graph, which has more
vertices. We summarize the results of this section in the following
theorem.  

\begin{theorem}\label{thm:swongbounds} 
  Let $\mathcal{M}^3_n$ be the set of 3-connected planar maps with $n$
  vertices and $\mathcal{S}(M)$ denote the set of Schnyder woods of
  $M\in\mathcal{M}^3_n$. Then,\[3.209^n\leq \max_{M\in\mathcal{M}^3_n}|\mathcal{S}(M)|\leq 8^{n}.\] 
\end{theorem}  
The example used for the lower bound is the square grid $G_{k,\ell}$.

\begin{theorem} For $k,\ell$ big enough the number of Schnyder woods
  of the augmented grid 
  $G^*_{k,\ell}$  is asymptotically
  $|\mathcal{S}(G_{k,\ell}^*)|\approx3.209^{k\ell}.$
\end{theorem}
\proof The graph induced by the non-rigid edges in the primal dual
completion $\widetilde{G}^*_{k,\ell}$ of $G^*_{k,\ell}$ is
$G_{2k-1,2\ell-1}-(2k-1,1)$. This is a 
square grid of roughly twice the size as the original and with the
lower left corner removed.
The rigid edges can be identified  using the fact that
$\alpha_S(v_\infty)=0$ and deleting them induces $\alpha_S'$ on 
$G_{2k-1,2\ell-1}-(2k-1,1)$. The new $\alpha_S'$ only differs from
$\alpha_S$ for vertices, which are incident to an outgoing rigid edge, and it
turns out, that $\alpha_S'(v)=d(v)-1$ for all primal or dual
vertices and  $\alpha_S'(v)=1$ for all edge vertices of
$G_{2k-1,2\ell-1}-(2k-1,1)$. Thus, a bijection between
$\alpha_S'$-orientations and perfect matchings of
$G_{2k-1,2k-\ell}-(2k-1,1)$ is established by identifying matching
edges with edges directed away from edge vertices. The closed form
expression for the number of perfect matchings of
$G_{2k-1,2k-\ell}-(2k-1,1)$ is known (see~\cite{kpw-tam-00}) to be
$$\prod_{i=1}^{k}\prod_{j=1}^{\ell}\left(4-2\cos\frac{\pi
    i}{k}-2\cos\frac{\pi j}{\ell}\right). $$ 

The number of perfect matchings of $G_{2k-1,2\ell-1}-(2k-1,1)$ is
sandwiched between that of $G_{2k-2,2\ell-2}$ and that of
$G_{2k,2\ell}$. Therefore the asymptotic behavior is the same and
in~\cite{lovaszplummer}, the limit of the number of perfect matchings of
$G_{2k,2\ell}$, denoted as $\Phi(2k,2\ell)$, is calculated to be
\[\lim_{k,\ell\rightarrow\infty}\frac{\log\Phi(2k,2\ell)}{2k\cdot
  2\ell}=\frac{\log2} 
{2}+\frac{1}{4\pi^2}\int_0^\pi\int_0^\pi\log(\cos^2(x)+\cos^2(y))dx
dy\approx 0.29.\]
This implies that $G_{k,\ell}^*$ has asymptotically
$e^{4\cdot0.29\cdot k\ell}\approx3.209^{k\ell}$ Schnyder woods.\qed

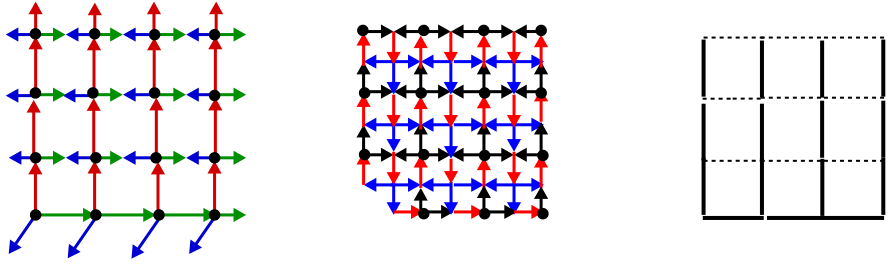
\PsFigCap{40}{swgridtree}{A Schnyder wood on  $G^*_{4,4}$, the
  reduced primal dual completion $G_{7,7}-(7,1)$ with the
  corresponding orientation and the associated spanning tree.} 

\ni\textbf{Remark.} In~\cite{temperley} Temperley discovered a
bijection between spanning trees of $G_{k,\ell}$ and perfect matchings
of $G_{2k-1,2\ell-1}-(2k-1,1)$. Thus, Schnyder woods of
$G_{k,\ell}^*$ are in bijection with spanning trees of
$G_{k,\ell}$, see Figure~\ref{fig:swgridtree}. This bijection can be
read off directly from the 
Schnyder wood: the unidirected edges not incident to a special
vertex form exactly the related spanning tree. Encoding both,
the Schnyder woods and the spanning trees, as $\alpha$-orientations
also gives an immediate proof of this bijection.\medskip

We now turn to the proof of the upper bound stated in
Theorem~\ref{thm:swongbounds}. The proof uses the upper
bound for Schnyder woods on plane triangulations, see
Theorem~\ref{thm:swontbounds}. We define a
triangulation $T_M$ such that there is in injective 
mapping of the Schnyder woods of $M$ to the Schnyder woods of
$T_M$. The triangulation $T_M$ is obtained from $M$ by adding a vertex
$v_F$ to every face $F$ of $M$ with $|F|\geq4$, see
Figure~\ref{fig:subdivideface}. The generic structure of a bounded
face of a Schnyder wood is shown on the left in the  top row of
Figure~\ref{fig:subdivideface}, for a proof see~\cite{f-gga-04}. The
three edges, which do not lie on the boundary of the triangle, are
the \term{special edges} of $F$.

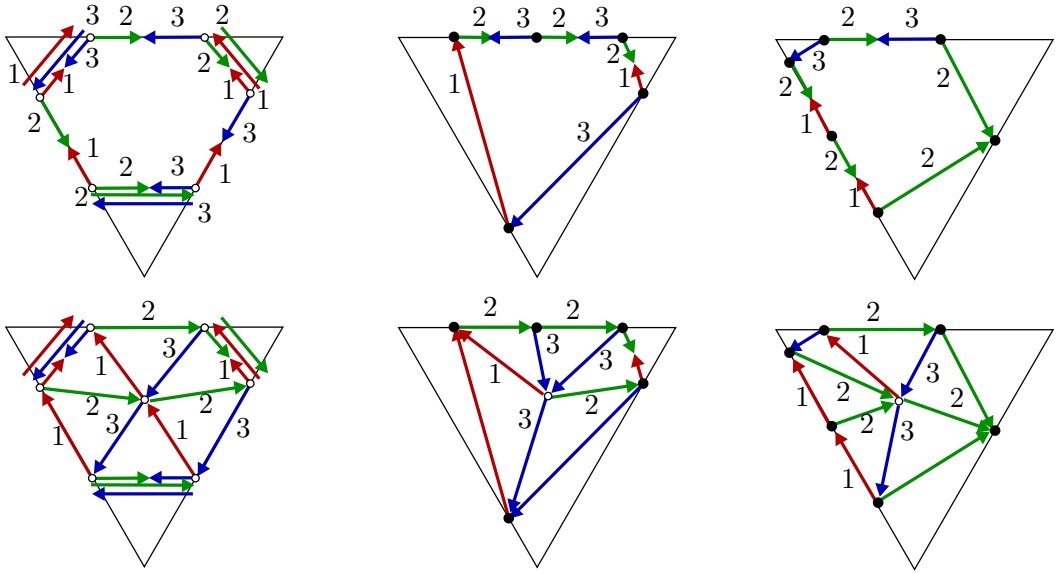
\PsFigCap{90}{subdivideface}{A Schnyder wood on a map $M$ induces a
  Schnyder wood on $T_M$. The three special edges of a face are those,
  which do not lie on the black triangle.} 

A vertex $v_F$ is adjacent to all the vertices of $F$. A Schnyder wood
of $M$ can be mapped to a Schnyder wood of $T_M$ using the generic 
structure of the bounded faces as shown in
Figure~\ref{fig:subdivideface}. The green-blue non-special edges of
$F$ become green unidirected. Their blue parts are substituted by
unidirected blue edges pointing from their original start-vertex
towards $v_F$. Similarly the blue-red non-special edges become blue
unidirected  and the red-green ones red unidirected. Three of the edges incident to $v_F$ are still undirected at
this point. They are directed away from $v_F$ and colored in
accordance with (W3).

 Let two different
Schnyder woods be given that have different directions or colors on an
edge $e$. That the map is injective can be verified by comparing the
edges on the boundary of the two triangles on which the edge $e$ lies
in $T_M$.

Thus, it suffices to bound the number of Schnyder woods of $T_M$. We
do this by specializing Proposition~\ref{prop:generalupperbound}. We
denote the set of vertices of $T_M$ that correspond to faces of size 4
in $M$ by $F_4$ and its size by $f_4$ and similarly $F_{\geq5}$ and
$f_{\geq5}$ are defined. Note that $I=F_4\cup F_{\geq5}$ is an
independent set and $T_M$ has a spanning tree in which all the
vertices from $I$ are leaves. Let $n_T$ denote the number of vertices
of $T_M$. Then, $T_M$ has at most 
\begin{equation} 2^{3n_T-6-n_T}\cdot\prod_{v\in I}
\left(\frac{1}{2^{d(v)-1}}\bin{d(v)}{3}\right) \leq
4^{n+f_4+f_{\geq5}}\cdot\left(\frac{1}{2}\right)^{f_4}
\cdot\left(\frac{5}{8}\right)^{f_{\geq5}}=
4^{n}\cdot2^{f_4}\cdot\left(\frac{5}{2}\right)^{f_{\geq5}}\end{equation}
Schnyder woods. Note that $n+f_4+f_{\geq5}+f_4+2f_{\geq5}\leq
m+f_4+2f_{\geq5}\leq 3n-6$ which implies that
$f_4+\frac{3}{2}f_{\geq5}\leq n$. Maximizing 
equation~(\arabic{equation}) under this condition yields that the
maximum $8^n$ is attained when $f_4=n$. Thus $M$ has no more than $8^n$
Schnyder woods.\qed

The proof of the lower bound $3.209^n$ involves the result about the
number of perfect matchings of the square grid. This result makes use
of non-combinatorial methods. Therefore, we complement this bound with
a result for another graph family, which uses a straight-forward analysis,
but still yields that these graphs have more Schnyder woods than the
triangular grid, see Section~\ref{sec:swont}. 

The graph we consider is the \term{filled hexagonal grid}
$H_{k,\ell}$, see Figure~\ref{fig:hexgridsym2}. Neglecting boundary
effects he hexagonal grid has twice as many  vertices as hexagons.
This can be seen by associating with every hexagon the 
vertices of its northwestern edge. Thus, neglecting boundary effects,
the filled hexagonal grid has five vertices per hexagon. The boundary
effects will not hurt our analysis because $H_{k,\ell}$ has only
$2(k+\ell)$ boundary vertices but
$5\cdot k\ell+2(k+l)$ vertices in total.

\begin{proposition} For $k,\ell$ big enough the filled hexagonal grid
  $H_{k,\ell}$ has
  \[2.63^{n}\leq|\mathcal{S}(H_{k,\ell})|\leq6.07^n.\]
\end{proposition}

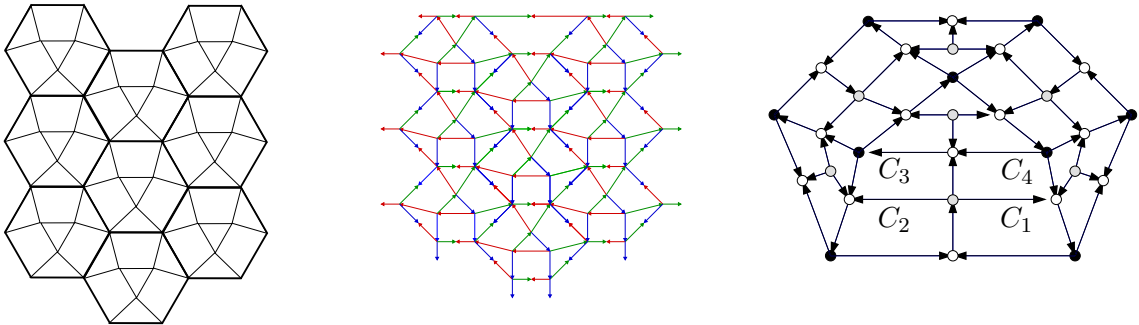
\PsFigCap{25}{hexgridsym2}{The filled hexagonal grid $H_{3,3}$, a
  Schnyder wood on this grid and the
  primal-dual suspension of a hexagonal building block of $H_{k,\ell}$. Primal
  vertices are black, face vertices blue and edge vertices green.}

\proof  We count how many different orientations we can have
on a filled hexagon. We do this using the bijection from
Theorem~\ref{thm:swong<->or}. The right part of
Figure~\ref{fig:hexgridsym2} shows a feasible
$\alpha_S$-orientation of a filled
hexagon. Note that this orientation is feasible on the boundary when
we glue together the filled hexagons to a grid $H_{k,\ell}$ and add
a triangle of three special  vertices around the grid. We flip only boundary
edges of a hexagon which belong to a 4-face in this hexagon. 
As these edges belong to a triangle in the hexagon on their
other side, the cycle flips in any two filled hexagons can be
performed independently. 

Let us now count how many orientations a filled hexagon
admits, see the right part of Figure~\ref{fig:hexgridsym2} for the
definition of the cycles $C_1, C_2, C_3$ and $C_4$. If the 6-cycle
induced by the central triangle is directed as shown in the rightmost
part of Figure~\ref{fig:hexgridsym2} then we can flip either $C_1$ or
$C_2$ and if $C_2$ is flipped, $C_3$ can be flipped as well. This yields $4^3$
orientations, as the situation is the same at the other two 4-faces of
the hexagon. If the 6-cycle is flipped the same calculation can be
done with $C_3$ replaced by $C_4$. This makes a total of
$2\cdot 4^3=128$ orientations per filled hexagon. That is, there are
at least $128^{k\cdot\ell}\geq 2.639^{5\cdot k\cdot\ell}$ orientations
of  $H_{k,\ell}$. 

We start the proof of the upper bound by collecting some statistics
about $H_{k,\ell}$.  As mentioned above, $H_{k,\ell}$ has $n=5\cdot
k\cdot\ell$ interior vertices, $12\cdot 
k\ell$ edges and $7 \cdot k\ell$ faces. Thus, the primal-dual completion has
$48 \cdot k\ell$ edges. There is no choice for the orientation of
the edges incident to the $3\cdot4/7\cdot
f=12\cdot k\ell$ face vertices of triangles. We can choose a
spanning tree $T$ on the remaining 
$5\cdot k\ell+12\cdot k\ell+3\cdot k\ell$ vertices such that all
face vertices are leafs and proceed as in the proof of
Proposition~\ref{prop:generalupperbound}, but using that we know the
number of edges exactly. Since in the independent set of
the remaining face vertices all of them have degree 4 and required
out-degree 3, they contribute a factor of $1/2$ each. Thus, there are at most
$ 2^{(48-12-20)k\cdot\ell}\cdot2^{-3\cdot
  k\ell}=2^{13\cdot k\ell}\leq 6.07^n$ Schnyder woods on
$H_{k,\ell}$.\qed 


\section{Counting 2-Orientations}\label{sec:2-or}

%
Felsner et al.~\cite{FHKO07} present a theory of 2-orientations of plane
quadrangulations, which shows many similarities with the theories of
Schnyder woods for triangulations. A quadrangulation is a planar map
such that all faces have cardinality four. A \term{2-orientation of a 
quadrangulation} $Q$ is an orientation of the edges such that all
vertices but two non-adjacent ones on the outer face have out-degree~2.

In~\cite{ffnops-bij-07} it is shown that 2-orientations on
quadrangulations with $n$ inner quadrangles are counted by the 
Baxter-number $B_{n+1}$. Hence asymptotically there are about 
$8^n$ 2-orientations on quadrangulations with $n$ vertices. Tutte
gave an explicit formula for rooted quadrangulations. A bijective proof
of Tutte's formula is contained in the thesis of Fusy~\cite{f-cpmaa-07}.
The formula implies that asymptotically there are about $6.75^n$
quadrangulations on $n$ vertices.  The two results together yield that a
quadrangulation with $n$ vertices has on average about $1.19^n$ 
2-orientations.

 \PsFigCap{55}{torus-or}
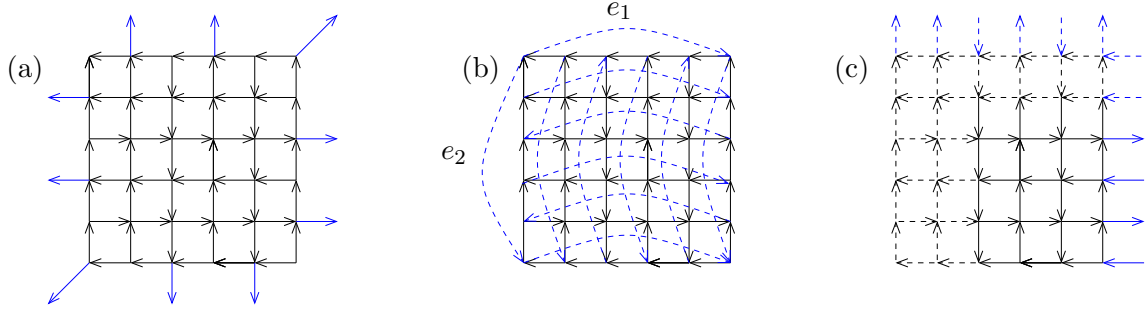
 {A 2-orientation of $G^\Box_{6,6}$, the corresponding Eulerian
   orientation $X$ of $G^T_{6,6}$ and an alternating orientation of
   $G^T_{4,4}$ that can be extended to $X$.}

We now give a lower bound for the number of 2-orientations of
$G^\Box_{k,\ell}$.   The proof method via transfer matrices and
eigenvalue estimates comes
from Calkin and Wilf~\cite{calkinwilf}. There it is used for
asymptotic enumeration of independent sets of the grid graph. Let
$\mathcal{Z}(Q)$ denote the set of all 2-orientations of a 
 quadrangulation $Q$, with fixed sinks.
 \begin{proposition}\label{prop:2-orlow}
  For $k,\ell$ big enough $G^\Box_{k,\ell}$ has
  \[1.537^{k\ell}\leq
  |\mathcal{Z}(G^\Box_{k,\ell})|\leq(8\cdot\sqrt{3}/9)^{k\ell}\leq
  1.5397^{k\ell}.\] 
\end{proposition}
\proof We consider 2-orientations of $G^\Box_{k,\ell}$ with sinks $(1,1)$ and
$v_\infty$. These 2-orientations induce Eulerian orientations
of $G_{k,\ell}^T$. The wrap-around edges inherit the direction of the
respective edges incident to $v_\infty$ (see
Section~\ref{sec:grids}) and $e_1,\;e_2$ are directed 
away from $(1,1)$. Therefore $G^\Box_{k,\ell}$ has at most as
many 2-orientations as $G_{k,\ell}^T$ has Eulerian orientations, which
implies the claimed upper bound.

Conversely a Eulerian orientation of $G_{k,\ell}^T$ in which the
wrap-around edges have these prescribed orientations induces a
2-orientation of $G^\Box_{k,\ell}$. Such Eulerian orientations are
called almost alternating orientations in the sequel,
see Figure~\ref{fig:torus-or}~(b).

Proving a lower bound for the number of almost alternating Eulerian
orientations yields a lower bound for the number of 2-orientations of
$G^\Box_{k,\ell}$.  

For the sake of simplicity we will work with alternating orientations
of $G^T_{k-2,\ell-2}$ instead of almost alternating ones of
$G_{k,\ell}^T$. In these Eulerian orientations the wrap-around edges
are directed alternatingly up and down respectively left and right,
see Figure~\ref{fig:torus-or}~(c). It is easy to see that this is a
lower bound for the number of almost alternating orientations of
$G_{k,\ell}^T$. Since we are interested in an  asymptotic lower bound
there is no difference in counting alternating orientations of
$G_{k-2,\ell-2}^T$ and  $G_{k,\ell}^T$ from our point of view and we
will continue working with alternating orientations of
$G_{k,\ell}^T$ to keep then notation simple.

Consider a vertex column $V_j^C$ of
$G_{k,\ell}^T$ and the edge columns $E^C_{j-1}$ and $E^C_j$. Let $X_1,\; X_2$
be orientations of $E^C_{j-1}$ respectively $E^C_j$. Let
$\delta(X_1,\;X_2)=1$ if and only if the edges induced by $V^C_j$ can
be oriented such that all the vertices of $V^C_j$ have out-degree~2. Let
$\delta_U(X_1,\;X_2)=1$ respectively $\delta_D(X_1,\;X_2)=1$ if and
only if $\delta(X_1,\;X_2)=1$ and the 
wrap-around edge induced by $V_j$ is directed upwards respectively
downwards. Note that
\[
\delta_U(X_1,\;X_1)=1=\delta_D(X_1,\;X_1)\] and 
\[
\delta_U(X_1,\;X_2)=1\Longleftrightarrow \delta_D(X_2,\;X_1)=1.\]

We define two transfer matrices $T_U(2k)$ and $T_D(2k)$. These are
square 0-1-matrices with the rows and columns indexed by the
$\bin{2k}{k}$ orientations of an edge column of size $2k$, that have
$k$ edges directed to the right. The transfer matrices are defined by
$(T_U(2k))_{X_1,X_2}=\delta_U(X_1,\;X_2)$ and
$(T_D(2k))_{X_1,X_2}=\delta_D(X_1,\;X_2)$. Hence $T_U(2k)=T_D(2k)^T$
and $T_{2k}=T_U(2k)\cdot T_D(2k)$ is a real symmetric non-negative
matrix with positive diagonal entries. From the combinatorial
interpretation it can be seen that $T_{2k}$ is primitive, that is
there is an integer $\ell\geq1$ such that all entries of $T_{2k}^\ell$
are positive and thus the Perron-Frobenius Theorem can be
applied. Hence, $T_{2k}$ has a unique eigenvalue $\Lambda_{2k}$ with
largest absolute value, its eigenspace is 1-dimensional and the
corresponding eigenvector is positive.

 Let $X_A$ be one of the two edge column
orientations that have alternating edge directions and $e_A$ the
vector of dimension $\bin{2k}{k}$ that has all entries 0 but the one
that stands for $X_A$, which is 1.

The number $c_A(2k,2\ell)$ of alternating orientations of $G^T_{2k,2\ell}$ is
$\left({T_{2k}^\ell}\right)_{X_A,X_A}= \langle e_A,T_{2k}^\ell 
e_A\rangle$. Since the eigenvector belonging to $\Lambda_k$ is
positive it is not orthogonal to any column of $T_{2k}$ and we obtain 
\[
\lim_{\ell\rightarrow\infty}c_A(2k,2\ell)^{1/\ell}=
\lim_{\ell\rightarrow\infty}\left(\left(T_{2k}^\ell\right)_{X_A,X_A}\right)^{1/\ell}
=\Lambda_{2k},
\]
where the last equality is justified by an argument known as the power
method.

It follows from~\cite{l-treosi-67} that the limit
$\lim_{k\rightarrow\infty}\Lambda_{2k}^{1/k}$ exists, but for the sake of
completeness we provide an argument from~\cite{calkinwilf}. We use
that $\Lambda_{2k}^p\geq \langle v,T_{2k}^p v\rangle/\langle 
v,v\rangle$ for any vector $v$ and that $\langle e_A,T_{2k}^{p}
e_A\rangle=\langle e_A,T_{2p}^k e_A\rangle$  
since both expressions count the number of alternating orientations of
$G^T_{2k,2p}$. 
 \[ \left(\Lambda_{2k}^{1/k}\right)^p 
 =\Big(\Lambda_{2k}^p\Big)^{1/k}
 \geq \Big(\langle e_A,T_{2k}^p e_A\rangle\Big)^{1/k}
 = \left(\langle e_A,T_{2p}^{k} e_A\rangle\right)^{1/k}
 \]
Taking limits with respect to $k$ on both sides yields
$$
 \left(\liminf_{k\rightarrow\infty}\Lambda_{2k}^{1/k}\right)^p
  \geq \liminf_{k\rightarrow\infty}\left(\langle e_A,T_{2p}^{k} e_A\rangle\right)^{1/k}
 = \Lambda_{2p}
$$
which implies 
$\displaystyle
\liminf_{k\rightarrow\infty}\Lambda_{2k}^{1/k}\geq\
 \limsup_{p\rightarrow\infty}\Lambda_{2p}^{1/p}.
$
It follows that $\displaystyle\lim_{k\rightarrow\infty}\Lambda_{2k}^{1/k}$ exists. Similar
arguments as above yield the following.
\[\Lambda_{2k}^{p}\geq \frac{ \langle e_AT^q_{2k},T^p_{2k}T^q_{2k} e_A\rangle} 
{\langle  T^q_{2k}e_A,T^q_{2k}e_A\rangle}= 
 \frac{ \langle e_A,T_{2k}^{p+2q}e_A\rangle} {\langle
   e_A,T_{2k}^{2q}e_A\rangle} =
\frac{ \langle e_A,T^{k}_{2p+4q}e_A\rangle} {\langle
  e_A,T^{k}_{4q}e_A\rangle}. 
\]
Taking limits with respect to $k$ on both sides yields
\[
\lim_{k\rightarrow\infty}\Lambda_{2k}^{1/k}
\geq\left(\frac{\Lambda_{4q+2p}}{\Lambda_{4q}}\right)^{1/p}   
\]
We are interested in
$\displaystyle\lim_{k\rightarrow\infty}
\lim_{\ell\rightarrow\infty}c_A(2k,2\ell)^{1/4k\ell}=
\lim_{k\rightarrow\infty}\Lambda_{2k}^{1/4k}$ since $4k\ell$ is the
number of vertices of $G^T_{2k,2\ell}$. Using a {\sl Mathematica}
program we have computed $\Lambda_{10}$ and 
$\Lambda_8$ with the result that 
\[\left(\frac{\Lambda_{10}}{\Lambda_{8}}\right)^{1/4}
\geq \left(\frac{2335.8714}{418.2717}\right)^{1/4} \geq1.537\]  
\qed

\ni\textbf{Remark.} 
We return to the correspondence between 2-orientations of
$G^\Box_{k,\ell}$ and Eulerian
orientations of  $G^T_{k,\ell}$, that was mentioned at the beginning
of the last proof.
By the pigeon hole principle, there must
be a sequence of orientations $X_{k,\ell}$ of the wrap-around edges
that extends asymptotically to $(8\cdot\sqrt{3}/9)^{k\ell}$ Eulerian
orientations of $G^T_{k,\ell}$. This implies that for $k,\ell$
big enough there is an $\alpha_{k,\ell}$ on $G_{k,\ell}$ such that
there are $(8\cdot\sqrt{3}/9)^{k\ell}$ $\alpha_{k,\ell}$-orientations
of  $G_{k,\ell}$. This $\aa_{k,\ell}$ satisfies $\alpha_{k\ell}(v)=2$
for every inner vertex $v$ and  $\alpha_{k\ell}(w)\in\{0,1,2\}$ for
every boundary vertex $w$. We call $\alpha$-orientations of this type
\term{inner 2-orientations of the grid}.

We think that $G^\Box_{k,\ell}$
has asymptotically $(8\cdot\sqrt{3}/9)^{k\ell}$ 2-orientations. But we
were not able to show this, just like for the case of the
triangular grid, see the last remark of Section~\ref{sec:swont}. 


\begin{theorem}\label{thm:2orbounds} Let $\mathcal{Q}_n$ denote the
  set of all plane quadrangulations with $n$  vertices and
  $\mathcal{Z}(Q)$ the set of 2-orientations of
  $Q\in\mathcal{Q}_n$. Then, for $n$  big enough 
  \[1.53^n\leq \max_{Q\in\mathcal{Q}_n}|\mathcal{Z}(Q)|\leq1.91^{n}.\]
\end{theorem}



\proof The lower bound is that from Proposition~\ref{prop:2-orlow}. An
upper bound of $2^n$ follows immediately from 
Lemma~\ref{lem:tree}. Note that we may assume that $Q$ does not have
vertices of degree 2, because their incident edges would be
rigid. Bonsma~\cite{b-phd-06,b-stwmligwmdt-07} shows that
triangle-free graphs of minimal degree 3 have a spanning tree $T$ with
more than $n/3$ leafs. As $Q$ is bipartite, $T$ has as a set $I$ of at
least $n/6$ leafs which is an independent set of $Q$. As in 
Proposition~\ref{prop:generalupperbound}, this yields that there are
at most $2^n\cdot (3/4)^{n/6}\leq 1.91^n$ 2-orientations of $Q$.\qed



\section{Counting Bipolar Orientations}\label{sec:bipo}

%
We first give an overview of the definitions and facts about bipolar
orientations that we need in this section. A good starting point for
further reading about bipolar orientations is~\cite{fraysseix95}. 

Let $G$ be a connected graph and $e=st$ a distinguished edge of 
$G$. An orientation $X$ of the edges of $G$ is an \term{$e$-bipolar
orientation} of $G$ if it is acyclic, $s$ is the only vertex without
incoming edges and $t$ is the only vertex without outgoing edges. We
call $s$ and $t$ the source respectively sink of $X$. There are many
equivalent definitions of bipolar orientations,
c.f.~\cite{fraysseix95}. 

Figure~\ref{bipolar} shows examples of plane bipolar orientations. 
Note that in this case it is enough to have  vertices $s$, $t$ on the outer face,
they need not be adjacent.  The following characterization of
plane bipolar orientations will be useful to keep some proofs in the
sequel simple.
\begin{proposition}\label{prop:bundles} 
An orientation $X$ of a planar map $M$ with two special
vertices $s$ and $t$ on the outer face is a bipolar orientation if and
only if it has the following two properties.
  \begin{description}
  \item{(1)} Every vertex other than the source $s$ and the sink $t$, has
    incoming as well as \\ outgoing edges.
  \item{(2)} There is no directed facial cycle. 
  \end{description}
Furthermore, the following stronger versions of the above properties
hold for every bipolar orientation.
  \begin{description}
  \item{(1')} At every vertex other than the source and the
  sink, the incoming and outgoing edges form two non-empty bundles of
  consecutive edges.
  \item{(2')} The boundary of every face has exactly one sink and one source,
  i.e.\ consists of two directed paths.
  \end{description}
\end{proposition} 
We omit the proof that properties $(1)$ and $(2)$ imply that $X$ is a
bipolar orientation. The proof that every bipolar orientation has
properties $(1')$ and $(2')$ (and thus properties $(1)$ and $(2)$ as
well) can be found in in~\cite{woods82} or~\cite{tamassia86}.

The fact that properties $(1')$ and $(2')$ characterize bipolar
orientations of planar maps yields a bijection
 between bipolar orientations of a
map $M$ and 2-orientations of the angular map, i.e., of the map
$\widehat{M}$ on the vertex set $V\cup \mathcal{F}$, where where
$\mathcal{F}$ is the set of faces of $M$, and edges $\{v,f\}$ for all
incident pairs with $v\in V$ and $F\in \mathcal{F}$.
This bijection was first described by Rosenstiehl~\cite{rosenstiehl83}.

\PsFigCap{55}{bipolar}
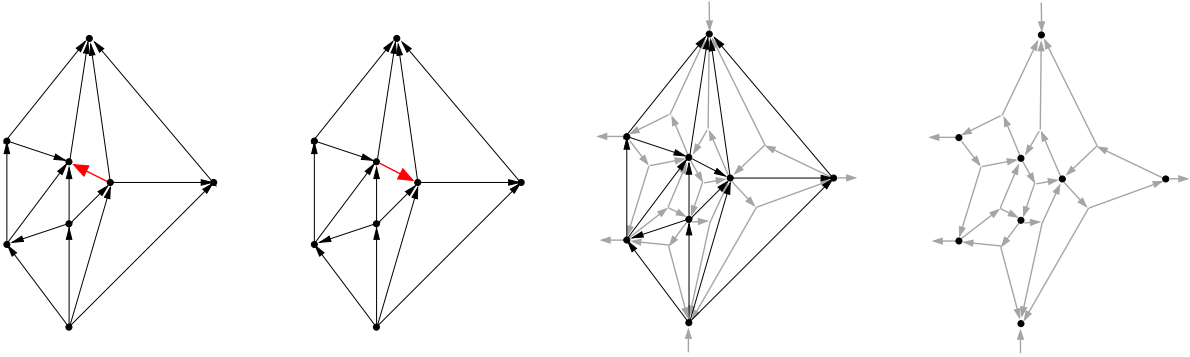
{Two bipolar orientations of the same graph with different out-degree
  sequences and the corresponding $\aa$-orientation of the angle graph
  (vertex for the unbounded face omitted).} 
Since bipolar orientations and 2-orientations of quadrangulations are
in bijection, we explain now how the results from
Theorems~\ref{thm:2orbounds} and~\ref{thm:bip} are related. A
triangulation with $n$ vertices has an angle graph with roughly $3n$
vertices. Hence the upper bound of $1.91^{3n}$, that
Theorem~\ref{thm:2orbounds} yields for the number of bipolar
orientation, is worse then the upper bound from
Theorem~\ref{thm:bip}. Conversely, every quadrangulation $Q$ with $n$ vertices is the
angle graph of the map obtained by connecting two vertices of
one of its partition classes by an edge, if they lie on a common
4-face of $Q$. This might yield a multi-graph if $Q$ has degree~2
vertices. But we may neglect this, since parallel edges must have the
same direction in every bipolar orientation. One of the partition
classes of $Q$ has size at most $n/2$, and thus  the upper bound from
Theorem~\ref{thm:bip} yields that $Q$ has at most $3.97^{n/2}$
2-orientations, which is worse than the bound from
Theorem~\ref{thm:2orbounds}. The example of the $G^\Box_{k,\ell}$,
which has $1.53^{n}$ 2-orientations is the angle graph of a graph which
has roughly $n/2=:n'$ vertices. Therefore this
yields only an example with $1.53^{2n'}$ bipolar orientations, which is
far away from the bound given in Theorem~\ref{thm:bip}. Conversely,
the triangular grid $T_{k,\ell}$, which has at least $2.91^n$ bipolar
orientations has an angle graph with roughly $3n=n'$ vertices. This
yields a quadrangulation with $2.91^{n'/3}$ 2-orientations, which is
worse than the bound for the number of 2-orientations that we obtained
in Proposition~\ref{prop:2-orlow}. 


\subsection{Counting Bipolar Orientations on the Grid}\label{sec:grid}


We now turn to analyzing the number of bipolar orientations of
$G_{k,\ell}$, with source $(1,1)$ and sink $(k,\ell)$ if $k$ is odd
and sink $(k,1)$ if $k$ is even.  For the proof of the Theorem, we
need sparse sequences. A \term{sparse sequence} is a 
0-1-sequences without consecutive 1s and it is well known that there
are $F_{n+2}$ such sequences of length $n$, where $F_{n+2}$ denotes
the $(n+2)$th Fibonacci number.

\begin{theorem}\label{thm:bipgrid}  Let and $\mathcal{B}(G_{k,\ell})$
  denote the set of bipolar orientations of   $G_{k,\ell}$. For
  $k,\ell$ big enough the number   of bipolar orientations of the grid
  $G_{k,\ell}$ is bounded by
  \[2.18^{k\ell}\leq|\mathcal{B}(G_{k,\ell})|\leq 2.619^{k\ell}.\]
\end{theorem}
\proof We first prove the lower bound with an argument using
directed cycles in a canonical orientation, as in
Proposition~\ref{prop:easylow}. To make this tool applicable we 
consider 2-orientations of the angle graph
$\hat{G}_{k,\ell}$ of $G_{k,\ell}$. Figure~\ref{fig:bipolargrid3}
shows the angle graph $\hat{G}_{4,5}$. The graph $\hat{G}_{k,\ell}$ has
$2k\ell-3(k+\ell)+4$ squares.  All dotted edges are rigid,
just like the four edges which are adjacent to a degree~2
vertex. Therefore, we may neglect all these edges in the rest of the proof.

The independent set $I$ of directed cycles in the canonical
orientation is 
marked by dots and includes approximately half of all squares. The
set $I'$ consists of all squares that are not in $I$. Members of $I'$
can be flipped if either the two cycles of $I$ 
above it, or the two cycles of $I$ below it are flipped, that is in 2
out of 16 cases (1 out of 4 for boundary squares). Roughly half of all
squares are in $I'$.

\PsFigCap{60}{bipolargrid3}
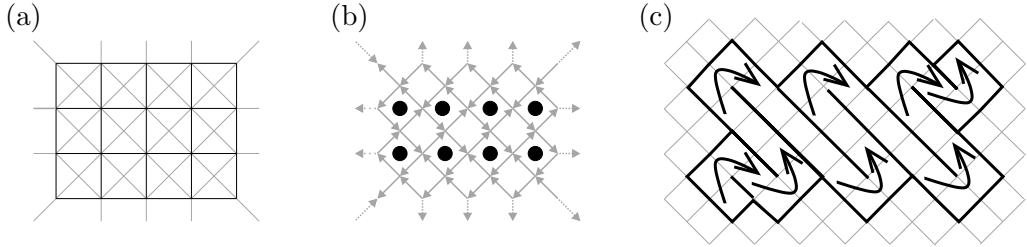
{ The grid $G_{4,5}$ is shown with its angle graph in blue; a
  canonical 2-orientation on $\hat{G}_{4,5}$ where red edges all
  connect to an additional vertex $v_\infty$ and the dots mark an
  independent set of directed cycles; the central part of
  $\hat{G}_{6,9}$ and the traversal used in the proof of the upper bound.}  
Thus, there are at least $2^{|I|+|I'|/8}$ bipolar orientations of
$G_{k,\ell}$, which leads to an asymptotic lower bound of
$2^{9k\ell/8}\approx 2.18^{k\ell}$.

For the proof of the upper bound we use a bijection discovered by
Lieb~\cite{l-treosi-67}. The bijection relates \term{face 3-colorings} where
no two squares sharing an edge have the same color and
inner 2-orientations of the square grid as shown in
Figure~\ref{fig:rgbrules2}~(a). Figure~\ref{fig:rgbrules2}~(b) shows
the face 3-coloring corresponding to the canonical 2-orientation, that
we used for the proof of the lower bound.

\PsFigCap{45}{rgbrules2}
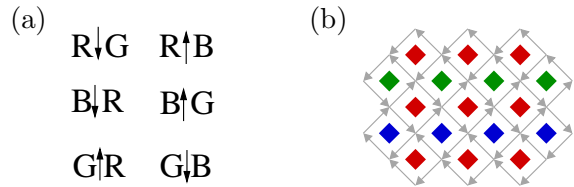
{Lieb's bijection between inner 2-orientations and face 3-colorings on
the grid.}  
Here we use this relation on $\hat{G}_{k,\ell}$.
We prove an upper bound for the number of face 3-colorings of
$\hat{G}_{k,\ell}$. The right part of Figure~\ref{fig:bipolargrid3}
shows the \emph{central part} of $\hat{G}_{k,\ell}$ bounded by a thick polygon.
We will encode the 3-coloring on the faces of the central part of
$\hat{G}_{k,\ell}$ as a sparse 
sequence $a$, where $a_i$ represents the $i$th square on the path $P$
indicated by the arrows in the figure.

The set $\mathcal{D}$ of faces,
which are not in the central part, has less than
$3^{|\mathcal{D}|}$ 3-colorings. In the encoding described next, the
code for the $i$th face of the 
path $P$ depends only on faces in $\mathcal{D}$ and faces of $P$ with index
smaller than $i$. Figure~\ref{fig:rgbencoding2} shows how the color of
the highlighted face is encoded by a 0 or a 1. The arrows indicate the direction
in which we traverse the central part of the graph. There are three
cases, one for a face where the path makes 
no turn and two for the two different types of turn faces. The
variables $X,Y,Z$ represent an arbitrary permutation of $R,G,B$. 

As for the decoding, it is clear from the figure that the faces marked
with an $X$ or $Y$ 
plus the 0-1 encoding uniquely determine the color of the face in
question. Thus the encoding is injective. It remains to show that
there cannot be consecutive 1s in this sequence. This follows from the
observation, that writing a 1 means, that the two faces that will be
used for the encoding of the next face on the path have different
colors. Thus, this face will be encoded by a 0.

\PsFigCap{55}{rgbencoding2}
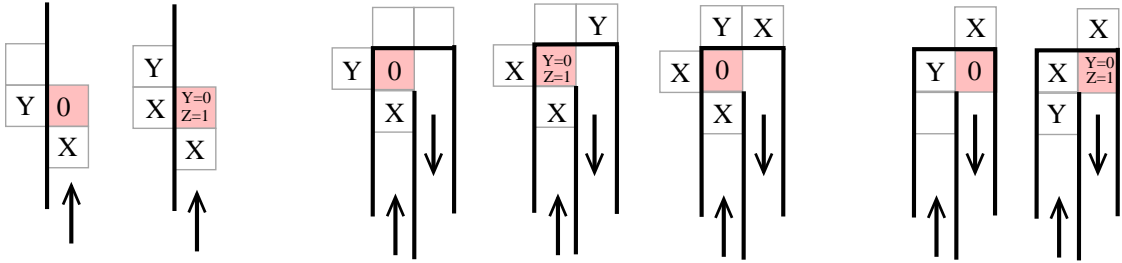
{Encoding a 3-coloring by a sparse 0-1$-$sequence. On the left the
  encoding for a square where the path makes no turn, in the center and
  right for the two different kind of turn faces.}
We bound the number of such encodings from above. The set
$\mathcal{D}$ can be covered by at most four 
horizontal plus four vertical rows of faces, thus $|\mathcal{D}|\leq
4(k+\ell)$. The length of the path is bounded by the number of bounded
faces of $\hat{G}_{k,\ell}$, which is less than $2k\ell$. Therefore,
there are at most \[3^{4(k+\ell)}\cdot F_{2k\ell+2}\] such
encodings. Using the asymptotics for the Fibonacci numbers this
implies, that there are less than $2.619^{k\ell}$ such encodings for
$k,\ell$ big enough.\qed

Lieb's analysis of the number of Eulerian
orientations of $T_{k\ell}$ is of interest in this case as well. 
It allows to improve the upper bound for grids with side lengths ratio
one to two.

\begin{proposition} For $k$ big enough the number of bipolar
  orientations of the grid $G_{k,2k}$ is bounded by
  \[2.18^{2k^2}\leq|\mathcal{B}(G_{k,2k})|\leq 2.38^{2k^2}.\]     
\end{proposition}
{\def\SetFigFont#1#2#3#4#5{%
  \fontsize{#1}{#2pt}%
  \fontfamily{#3}\fontseries{#4}\fontshape{#5}%
  \selectfont}%
\PsFigCap{53}{torus2aztek}
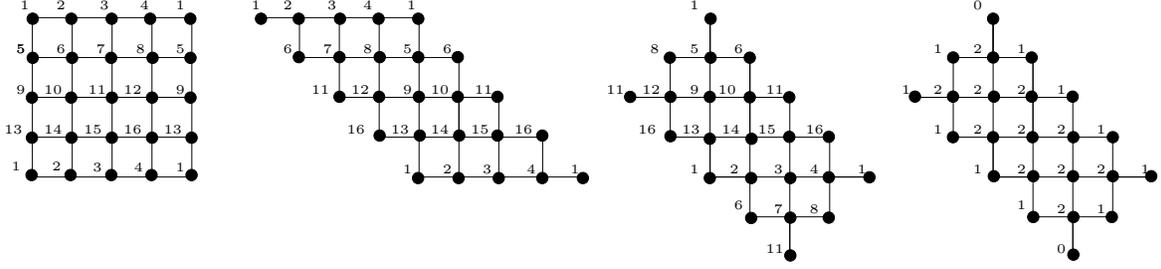
{How to obtain the tilted grid $\hat{G}_{5,3}$ from $G^T_{4,4}$ with
  two cuts. The numbers in the first three drawings are the vertex
  labels, in the last one they indicate $\hat{\alpha}$.} 
} 

\proof By  $\hat{G}_{k,\ell}'$ we denote the graph obtained from
$\hat{G}_{k,\ell}$ by deleting $v_\infty$ and all incident edges,
which are shown as dotted edges in Figure~\ref{fig:bipolargrid3}.
Figure~\ref{fig:torus2aztek} shows how to cut $G^T_{4,4}$ in 
two steps such that the grid looks 
like $\hat{G}_{3,5}'$ (if we do not identify vertices). The last
drawing shows, that every $\hat{\alpha}$-orientation of 
$\hat{G}_{3,5}'$ yields a Eulerian orientation of $G^T_{4,4}$ when we do
the appropriate identifications. In
general, this approach yields an injection from the bipolar
orientations of $G_{k+1,2k+1}$ to the Eulerian orientations of
$G^T_{2k,2k}$. As Lieb~\cite{l-treosi-67} has shown that $G^T_{2k,2k}$ has
asymptotically $(8\cdot\sqrt{3}/9)^{4k^2}$ Eulerian orientations,
this yields an upper bound of $(64/27)^{2k^2}$ for the number of
Eulerian orientations of $G_{k+1,2k+1}$. Every bipolar orientation of
$G_{k,2k}$ can be complemented to a bipolar orientation of
$G_{k+1,2k+1}$, thus $G_{k,2k}$ has at most as many bipolar
orientation as $G_{k+1,2k+1}$. The lower bound follows from the more
general claim of Theorem~\ref{thm:bipgrid}. \qed 

\ni\textbf{Remark.} The same problems as described in the  closing
remarks of Sections~\ref{sec:swont} and~\ref{sec:2-or} arise here when
trying to show that $G_{k,2k}$ actually has $(64/27)^{2k^2}$ bipolar
orientations by using Lieb's result for the torus.\medskip


\subsection{Counting Bipolar Orientations On Planar Maps}\label{bipolar}

Note that adding edges to the faces of size at least 4 of a planar map
$M$ can only increase the number of bipolar orientations by
Proposition~\ref{prop:bundles}. Thus we can restrict our considerations to
plane inner triangulations in this section. 
\begin{theorem}\label{thm:bip} Let $\MM_n$  denote the
  set of all planar maps with $n$ vertices and $\mathcal{B}(M)$ the
  set of all bipolar orientations of $M\in\MM_n$. Then, for $n$  big
  enough
  \[2.91^n\leq \max_{M\in\MM_n}|\mathcal{B}(M)|\leq3.97^{n}.\]
\end{theorem}
For the proof we need a couple of facts about Fibonacci numbers,
which are summarized in the following lemmas. The Fibonacci numbers
are the integer series defined by the recursion
$$
  F_1=1,\;F_2=1,\;F_n=F_{n-1}+F_{n-2}\mbox{ for }n\geq3.
$$
Define $F_0=0$ and let $\phi=\frac{1+\sqrt{5}}{2}$ be the Golden Ratio.

\begin{lemma}\label{lem:fib1} The Fibonacci numbers have the following
  properties.  
  \begin{itemize}
    \item $F_n=\frac{\phi^n-(1-\phi)^n}{\sqrt{5}}$ 
    \item $\lim_{n\rightarrow\infty}F_n=\frac{\phi^n}{\sqrt{5}}$ 
    \item $\sum_{i=0}^nF_iF_{n-i}=
      \frac{1}{5}\left(n(F_{n+1}+F_{n-1})-F_n\right)$  
  \end{itemize}
\end{lemma}
The first two are standard results from the vast theory of Fibonacci
numbers. The last formula is attributed to Shiwalkar and Deshpande 
in~\cite[A001629]{onlineseq}. The next lemma summarizes facts  about
sparse sequences. 

\begin{lemma}\label{lem:fib2} The number of sparse sequences of
  length $n$ is $F_{n+2}$. Let $r_n(i)$ be the number of sparse
  sequences of length $n$ whose $i$th entry is 1. Then,
\begin{itemize}
\item $r_n(i)=F_i\cdot F_{n+1-i}$
\item
  $\displaystyle{\sum_{i=1}^nr_n(i)=
    \frac{1}{5}\left(2(n+1)F_n+nF_{n+1}\right)}$  
\item
  $\displaystyle{\lim_{n\rightarrow\infty}
    \frac{\sum_{i=1}^nr_n(i)}{nF_{n+2}}= 
    \frac{1}{\sqrt{5}\phi}} \approx 0.2764$
\end{itemize}
\end{lemma}
The first identity follows from a construction of sparse
sequences of length $n$ from sparse sequences of length $n-1$ plus the
string ``0'' and sparse sequences of length $n-2$ plus ``01''. The
second and third identity then follow using the facts from
Lemma~\ref{lem:fib1}. 

Before proving the  Theorem~\ref{thm:bip} we give two results for the number of
bipolar orientations of special classes of planar maps. 
\begin{proposition}\label{prop:stack} A stacked triangulation with $n$
  vertices has $2^{n-3}$ bipolar orientations. 
\end{proposition}
\proof The  $K_4$ has two bipolar orientations for fixed source and
sink. We proceed by 
induction and assume, that a stacked triangulation with  $n$
vertices has $2^{n-3}$ bipolar orientations. Now let $T$ be a stacked
triangulation with $n+1$ vertices and $v$ a vertex of degree 3 in
$T$. Then, $T-v$ has $2^{n-3}$ bipolar orientations by induction. Now
stacking $v$ into $T$ again, there are exactly two ways to complete a
given bipolar orientation on $T-v$ without violating Properties $(1)$
or $(2)$ from Proposition~\ref{prop:bundles}. Thus,
there are $2^{(n+1)-3}$ bipolar orientations of $T$. \qed

\begin{proposition}\label{prop:outerbip} Let $\mathcal{O}_n$ be the
  set of all outerplanar maps  with $n$ vertices. Then,
  \[\max_{M\in\mathcal{O}_n}|\mathcal{B}(M)|= F_{n-1}\approx
  1.618^{n-1}.\]
\end{proposition}
\proof We show first that there are indeed outerplanar maps with
$F_{n-1}$ bipolar orientations. Let $T:=T_{2,\ell}$ be the
triangular grid with two rows. We consider bipolar orientations of
$T$ with source $(1,1)$ and sink $(2,\ell)$. In every such bipolar
orientation the boundary edges form two directed paths from
$(1,1)$ to $(2,\ell)$. We start by defining the standard bipolar
orientation $B_0$ of $G^*$, which is shown in
Figure~\ref{fig:outerplanar}. 

\PsFigCap{45}{outerplanar}
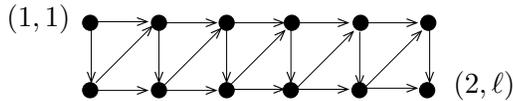
{The standard bipolar orientation on $T_{2,\ell}$.}
In $B_0$ the vertical inner edges are directed downwards and the
diagonal ones upwards. Now we encode any other orientation
of the inner edges by a  sequence $(a_i)_{i=1\ldots n'}$ of length
$n'=n-3$, where $a_i=1$ if the corresponding edge has the opposite
direction as in $B_0$ and $a_i=0$ otherwise. The entries come in the natural
left to right order in  $(a_i)_{i=1\ldots n'}$. We show that all
sparse sequences of length $n'$ produce bipolar 
orientations. In a sparse sequence there are no consecutive 1s, thus
out of the two inner edges incident to a vertex  at most one is
reversed with respect to $B_0$. This guarantees that there is no
directed facial 3-cycle. As all vertices have an incoming and an
outgoing outer edge the resulting orientation is bipolar, according to
Proposition~\ref{prop:bundles}. 

It remains to show that $F_{n-1}$ is an upper bound for the number of bipolar 
orientations of any outerplanar map $M$ with $n$ vertices. We may
assume that $M$ is a plane inner triangulation. The proof uses  
induction on the number of vertices and the claim is trivial for
$n=3$. Now, let $M$ have $n+1$ vertices and let $s$ be the source
vertex. If $M$ has a vertex $x\neq s,t$ of degree 2 with neighbors
$v,\; w$ then the direction
of the edge $\{v,w\}$ determines the directions of the edges
$\{x,v\}$ and $\{x,w\}$. Therefore, $M$ has at
most as many bipolar orientations as $M-x$, that is at most $F_{n-1}$
many. If all vertices but $s$ and $t$ have degree at least three, then $s$
and $t$ have degree 2 and the vertices of every inner edge of $M$ are
separated by $s$ and $t$ on the outer cycle. This is because the
interior of the boundary cycle on $n+1$ vertices is partitioned into $n-1$
triangles, and thus two of these triangles must share two edges with
the boundary, which yields two degree 2 vertices.

So $s$ is incident to
only two vertices $v$ and $w$, and we may assume that $v$ has degree 3
in $M$, that is the inner edge $e=\{v,w\}$ is the only
inner edge incident to $v$. Now, let $X$ be some bipolar orientation of
$M$ in which $e$ is directed from $v$ to $w$. Then, the orientation of
$M-s$ induced by 
$X$ is  a bipolar orientation with source $v$. For a bipolar
orientation $Y$ in which $e$ is oriented from $w$ to $v$, the
orientation of $M-s$ induced 
by $X$ is  a bipolar orientation with source $w$ and $v$ is a vertex of
degree 2 in $M-s$. This mapping is injective, and thus $M$ has at most
as many bipolar orientations as $M-s$ and $M-\{s,v\}$ together, that is
$F_{n-1}+F_{n-2}=F_{n}$. \qed

\ni\textbf{Remark.}  From the above proof it also follows, that
$T_{2,\ell}$ is the only outerplanar map on $2\ell$ vertices, which has
$F_{2\ell-1}$ bipolar orientations.\medskip

The example that gives the lower bound for the number of bipolar
orientations of planar maps is the triangular grid  $T_{k,k}$ with
source $(1,1)$ and sink $(k,k)$. 
\begin{proposition}\label{prop:bipolarlow2}
Let $T_{k,k}$ be the triangular grid and $k$ big enough. Then,
\[|B(T_{k,k})|\geq2.91^{n}.\]  
\end{proposition}
\proof We first claim that $T_{k,k}$ has at least $2.618^{k^2}$
bipolar orientations. To see this we glue together $k-1$ copies of
$T_{2,k}$. Every orientation of $T_{k,k}$ obtained in this way
corresponds to a concatenation of $k-1$ sparse sequences of length
$2k-3$, which we call an almost sparse sequence. We 
denote the set of all such sequences of length  $2k^2-5k+3$ by
$S$, the cardinality of $S$ is $F_{2k-1}^{k-1}$ which is bounded below by
$F_{2k^2-5k+3}\geq 2.618^{k^2}$ for $k$ big enough. That each $s\in S$
corresponds to a 
bipolar orientation of $T_{k,k}$ can be checked using
Proposition~\ref{prop:bundles}.  

The horizontal edge $e_{i,j}:=(i,j)\rightarrow(i,j+1)$ lies on the
boundary of two triangles for $2\leq i\leq k-1$. The other four edges
of these triangles are
\[\{(i,j),(i-1,j+1)\},\;\{(i,j+1),(i-1,j+1)\},\;\{(i,j),(i+1,j)\},\; 
\{(i,j+1),(i+1,j)\}.\] 
The crucial observation for improving the above bound is, that we can
reorient $e_{i,j}$ if and only if the entries belonging to these four
edges show one of the two patterns $10\ldots01$ or $01\ldots10$. 

We now choose $k-1$ sparse sequences of length $2k-3$ independently
uniformly at random 
and concatenate them to obtain a random almost sparse sequence $s\in S$.
It follows from the first identity from Lemma~\ref{lem:fib2} with
$n=2k-3$ and $i=2j-1$ 
that for $\{(i,j),(i,j+1)\}$ there are
$F_{2j-1}F_{2k-2j-1}F_{2k-1}^{k-2}$ sequences that have
$\{(i,j),(i-1,j+1)\}$ marked 1 out of the total $F_{2k-1}^{k-1}$
sequences. This is used to
calculate the probability that the entry for $\{(i,j),(i-1,j+1)\}$ is 1 as
\begin{eqnarray}\label{eqn2}
\lim_{k\rightarrow\infty}\frac{F_{2j-1}F_{2k-2j-1}F_{2k-1}^{k-2}}
{F_{2k-1}^{k-1}}&
=& \lim_{k\rightarrow\infty}\frac{1}{\sqrt{5}}
\cdot\phi^{2j-1}\cdot\phi^{2k-2j-1}\cdot\phi^{-2k+1}= \frac{1}{\sqrt{5}\phi}.
\end{eqnarray}
Taking the limit is only justifiable if $(2j-1)\rightarrow\infty$ and
$(2k-2j-1)\rightarrow\infty$ for 
$k\rightarrow\infty$. Therefore we introduce $\delta>0$ and denote the
set of horizontal edges with $\delta (k-1)\leq j\leq (1-\delta)(k-1)$ and 
$2\leq i\leq k-1$ by $E_\delta$. Taking the limit in equation (5) is
justified for all $e_{i,j}\in E_\delta$. The size of this set is
$|E_\delta|=(1-2\delta)(k-1)^2$ .

This works analogously for the edge $\{(i,j+1),(i+1,j)\}$ and the
events are independent. The pattern $01\ldots10$ has the same
probability in the limit and the patterns mutually exclude each
other. Thus, for every $\epsilon>0$ the probability that the edge
$\{(i,j),(i,j+1)\}\in E_\delta$ can be flipped is  
\[\PP[\Eins_{i,j}(s)=1]\geq \frac{2}{5\phi^2}-\epsilon\]
for $k$ big enough. 

The rest of the proof is to analyze how many of the flip-patterns we expect
for a sparse sequences $s$. Let  $Q(s)=\sum_{i,j}\Eins_{i,j}(s)$ be
a random variable 
counting the  number of flippable edges in $s$. We use Jensen's
inequality to estimate the number of orientations, which is at least 
\[|S| \cdot \EE_{s\in\mathcal{S}}[2^{Q(s)}] \geq 
F_{2k^2-5k+3}\cdot2^{\EE_{s\in\mathcal{S}}[{Q(s)}]}
\geq 2.618^{k^2}\cdot2^{
  (1-2\delta)(k-1)^2\left(\frac{2}{5\phi^2}-\epsilon\right)}\geq2.91^{k^2}
\]
for $k$ big enough. \qed

\ni\textbf{Remark.} The third identity from Lemma~\ref{lem:fib2} gives
the expected number of 1s that a random entry of a random sparse
sequence of length $n$ has. In equation~(\ref{eqn2}) we calculate the
expected number of 1s in a random sequence of length $n$ at a fixed
entry which is far enough away from the boundary. The fact that these
two values agree shows that the expected number of 1s at a fixed entry
of a sparse sequence does not depend strongly on the choice of the
entry.\medskip

The following relation is useful to upper bound the number of bipolar
orientations for plane inner triangulations. It has been
presented with a different proof in~\cite{m-bopg-06}. Let
$\mathcal{F}_b$ be the set  of bounded faces of $M$ and $\mathcal{B}$ the
set of bipolar orientations of $M$. Fix a bipolar orientation $B$.
The boundary of every triangle
$\Delta\in \mathcal{F}_b$ consists of a path of length two 
and a direct edge from the source to the sink of $\Delta$.
We say that $\Delta$ is a $+$ triangle of $B$ if looking along
the direct source-sink edge the triangle is on the left. Otherwise,
if the triangle is on the right of the edge we speak of a $-$
triangle, see Figure~\ref{fig:plusminus}.
We use this notation to define a mapping $G_B: \mathcal{F}_b \to \{-,+\}$.
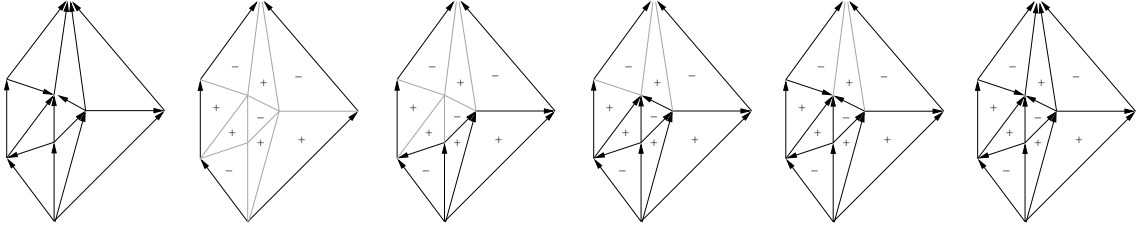
\PsFigCap{42}{plusminus}{A bipolar orientation, the corresponding
  $+/-$ encoding and an illustration of the decoding algorithm.}

\begin{theorem}\label{thm:plusminus} Let $M$ be a plane inner
  triangulation and $B$ a bipolar orientation of $M$.
  Given $G_B$, i.e., the signs of bounded faces, it is possible
  to recover $B$. In other words the function $B \to G_B$ 
  is injective from $\mathcal{B}(M) \to \{-,+\}^{|\mathcal{F}_b|}$.
\end{theorem}
\proof Given $G_B$ we construct $B$. We start by orienting all edges
on the boundary of the outer face such that $s$ and $t$ are the unique
source and sink of this face. We extend this partial orientation $Y$
with two rules. The {\em vertex rule} is applied to a vertex $v$, which
already has incoming and outgoing edges. It takes a bundle of
consecutive edges of $v$ which is bounded by two outgoing edges. 
It orients all the edges of the bundle such that they are outgoing at~$v$.
Note that these edge orientations are forced by Property $(1')$ of bipolar 
orientations in Proposition~\ref{prop:bundles}.
The {\em face rule} is applied to a facial triangle $\Delta$, which has
two oriented edges. The sign $G_B(\Delta)$ is used to deduce the 
orientation of the third edge.

Note that these two rules preserve the property that every vertex $v$,
which is incident  to an oriented edge in $Y$, can be reached from $s$
along an oriented path. In particular $v$ has an incoming edge.

Let $A_Y$ be the union of all faces, which have all boundary edges
oriented.  Initially, $A_Y$ consists of the outer face. Since $B$ is
acyclic the boundary of $A_Y$ is acyclic as well. Consequently, as
long as there are faces which do not belong to $A_Y$, there is a vertex
$v$ on the boundary of $A_Y$ which has two outgoing edges that belong
to the boundary of $A_Y$.  Either $v$ is a candidate to extend the
orientation using the vertex rule or there is a face incident to $v$
which becomes an element of $A_Y$ by applying the face rule to it. 

We have thus shown that the rules can be applied until $A_Y$ is the
whole plane, i.e., all edges are oriented. They have to be oriented as in $B$,
by construction. 
\qed

 The next theorem gives a necessary and sufficient condition for a
 vector in $\{-,+\}^{|\mathcal{F}_b|}$ to induce a bipolar
 orientation.  For the sake of simplicity we state it only for
 triangulations, but  the generalization to inner triangulations is
 straight forward. In order to obtain a more elegant formulation, we
 adopt the convention, that the unbounded face is signed $+$ if the
 bounded face adjacent to sink and source is signed $-$. Otherwise
 the unbounded face is signed $-$. Thus we work now with signings of
 the set $\mathcal{F}$ of all faces. We say that a $+$ triangle is the
 right knee of the vertex at which it has in incoming and an outgoing
 edge. Similarly a $-$ triangle is the left knee of exactly one of its
 vertices.

 \begin{theorem}\label{thm:plusminus2} Let $T$ be a triangulation,
   $x\in\{-,+\}^{|\mathcal{F}|}$, and $\mathcal{F}^-$ and 
   $\mathcal{F}^+$ the sets of faces that have negative
   respectively positive sign in $x$. Let $\widehat{M}(T)^+$ and
   $\widehat{M}(T)^-$ denote the subgraphs of the reduced angle 
   graph $\widehat{M}(T)-\{s,t\}$ induced by $V\cup \mathcal{F}^-$
    respectively $V\cup\mathcal{F}^+$. 

    Then, $x$ induces a bipolar orientation on $T$ if and only if both,
    $\widehat{M}(T)^+$ and $\widehat{M}(T)^-$, have a unique perfect
    matching.  
 \end{theorem} 

Theorem~\ref{thm:plusminus2} implies that every vertex of $T$ other than $s,t$  must be
adjacent to at least one $+$ triangle and one $-$ triangle. 

\begin{proposition} Let $T$ be a plane inner triangulation. Then, $T$
  has at most $3.97^n$ bipolar orientations.
\end{proposition}
\proof By Euler's  formula, there are $2^{f-1}\leq
2^{2n-2-f_\infty}=4^{n-1}\cdot2^{-f_\infty}$ many binary vectors 
of length $f-1$. By the
bijection from Theorem~\ref{thm:plusminus},
$4^{n-1}\cdot2^{-f_\infty}$ is also an upper bound for the number of
bipolar orientations of $T$.  

To squeeze the bound below $4^n$ we use the above observation that every
vertex of $T$ must be adjacent to at least one plus $+$ triangle and
at least one $-$ triangle. Thus, out of the $2^{d_v}$ possible 
$+/-$ vectors at a vertex $v$ at least two are not
feasible. Similarly, at an outer vertex $v\neq s,t$, there is
exactly one angle forming a knee at $v$. The sign of this angle
depends on which of the two oriented paths of the outer boundary $v$
lies on, but it is fixed either way. Thus, out of the $2^{d_v-1}$
possible sign patterns at $v$ at least one is not feasible. We
summarize, that at most a fraction of $(1-2^{1-d_v})$ of all sign
vectors is potentially feasible at every vertex but $s$ and $t$. We
denote the set $V\setminus\{s,t\}$ by $V'$ and its cardinality by
$n'$, i.e. $n'=n-2$. 

We apply Jensen's inequality, which says that for a convex
function $\varphi$ the inequality
$\varphi\left(\sum x_i/n\right)\leq\left(\sum
  \varphi(x_i)\right)/n$ holds and the inequality is reversed if
$\varphi$ is a concave function. As $\log x$ is concave we obtain
\[\log\left(\left(\prod_{v\in
      V'}(1-2^{1-d_v})\right)^{1/n'}\right)=\frac{1}{n'}\sum_{v\in
  V'}\log\left(1-2^{1-d_v}\right) \leq\log\left(\frac{1}{n'}\sum_{v\in
    V'}\left(1-2^{1-d_v}\right)\right).\]
By the monotonicity of the logarithm this implies
\[\prod_{v\in V'}(1-2^{1-d_v})\leq\left(\frac{1}{n'}\sum_{v\in
    V'}\left(1-2^{1-d_v}\right)\right)^{n'}.\]
The function $2^x$ is convex and we apply Jensen's inequality again
which yields
\[\frac{1}{n'}\sum_{v\in V'}\left(1-2^{1-d_v}\right)=1-2\left(\frac{1}{n'}\sum_{v\in
    V'}2^{-d_v}\right)\leq 1-2^{1+1/n'\cdot\sum_{v\in V'}-d_v}.\]
Since we deal with simple plane graphs $\sum_{v\in V'}d_v\leq
2(3n-6)=6n'$ and we conclude 
\[\prod_{v\in
  V'}(1-2^{1-d_v})\leq\left(1-2^{1-6}\right)^{n'}=\left(\frac{31}{32}\right)^{n'}.\]
Since we deal with plane graphs, $T$ can be partitioned into at
most four independent sets $I_k$, $k=1,\ldots,4$ by the Four Color Theorem. Thus, 
\[\left(\frac{31}{32}\right)^{n'}\geq\prod_{v\in
  V'}(1-2^{1-d_v})=\prod_{k=1}^4\prod_{v\in I_k}(1-2^{1-d_v})\] and for
at least one of the independent set it must hold that 
\[\left(\frac{31}{32}\right)^{n'/4}\geq\prod_{v\in I_k}(1-2^{1-d_v}).\]
We are ready to conclude that there are at most 
\[4^{n-1}\cdot2^{-f_\infty}\cdot\left(\frac{31}{32}\right)^{(n-2)/4}
<3.97^{n}\cdot2^{-f_\infty}\cdot\left(\frac{32}{31}\right)^{1/2}<3.97^{n}\]
bipolar orientations of $T$. 
\qed\medskip

\section{Complexity of Counting $\aa$-Orientations}

Given a planar map $M$ and some $\aa: V \rightarrow \NN$,
what is the complexity of computing the number of $\aa$-orientations of $M$?
In some instances this number can be computed efficiently, e.g. for
perfect matchings and spanning trees of general planar maps. In
Section~\ref{subsec:complexity} we show that counting  is
$\#{}P$-complete for other
$\aa$-orientations, and in
Section~\ref{subsec:perfmatch} we discuss how to adapt an existing fully
polynomial randomized approximation scheme.


\subsection{$\#{}P$-Completeness}
\label{subsec:complexity}



Recently, Creed~\cite{creed07} has shown the following theorem. As
already mentioned Eulerian orientations are $\alpha$-orientations 
and hence Theorem~\ref{thm:creed} says that counting $\aa$-orientations
is $\#{}P$-complete. Since we use the proof technique in the sequel,
we sketch the proof here. It uses a reduction from counting Eulerian
orientations, which has been proven to be $\#{}P$-complete by Mihail and Winkler 
in~\cite{mw-otnoeooag-96}. 

\begin{theorem}\label{thm:creed} It is $\#{}P$-complete to count
  Eulerian orientations of planar graphs. 
\end{theorem}
\proof
We aim to show that number Eulerian orientations of a graph $G$ can be
computed in polynomial time with the aid of polynomially many calls to
an oracle for the number of Eulerian orientations of a planar graph.

In order to count the Eulerian orientations of a graph $G$ with $n$
vertices a drawing of this graph in the plane with $\ell$ crossings is
produced. We may assume that no three edges cross in the same
point and that $\ell$ is of order~$O(n^2)$. 

From this drawing a family of graphs $G_i$ for
$i=0,\ldots,\ell$ is produced. In $G_i$ every crossing of two edges
$\{u,v\}$ and $\{x,y\}$ is replaced by the 
crossover box $H_i$ on $4k+1$ vertices, see
Figure~\ref{fig:planar-gadget}. For example in $G_0$ every crossing
is replaced by just one vertex. We call
the edges $\{v,w_k^u\},\;\{u,w_k^v\},\;\{x,w_k^x\},\;\{y,w_k^y\}$ the
\term{connection edges of $H_k$}.

Every
Eulerian orientation of $G$ induces a Eulerian orientation of $G_0$,
but there are Eulerian orientations of $G_0$, which do not come from a
Eulerian orientation of $G$. Given a Eulerian orientation of $G_0$ we call
the orientation of the edges incident 
to a vertex $w_0$, which replaces a crossing, valid if exactly one of
the edges $\{v,w_0\},\;\{u,w_0\}$ is 
directed away from $w_0$, and invalid otherwise. 

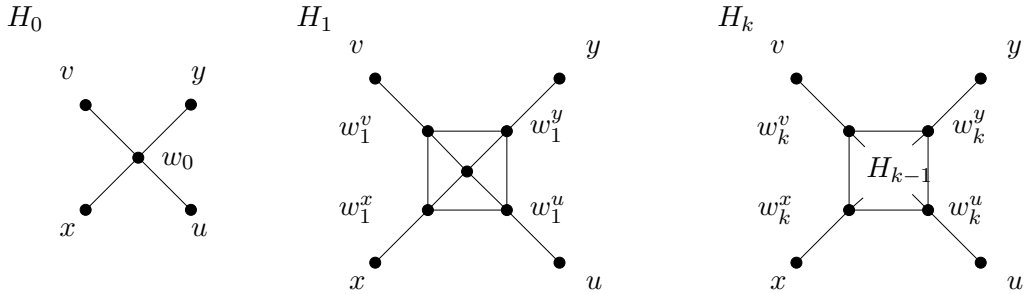
\PsFigCap{70}{planar-gadget}{The crossover boxes defined in~\cite{creed07}.}
We call a configuration of the connection edges of $H_k$ valid if
exactly one of the edges $\{v,w_k^v\},\;\{u,w_k^u\}$ is 
directed towards $\{u,v\}$ and exactly one of the edges
$\{x,w_k^x\},\;\{y,w_k^y\}$ is directed towards $\{x,y\}$. We call a
configuration of the connection edges of $H_k$ invalid if exactly two
of the connection edges are directed towards $\{u,v,x,y\}$, but it is
not a valid configuration. By $x_k$ respectively $y_k$ we denote the
number of ways that a valid 
respectively invalid configuration of the connections edges of $H_k$
can be extended to an orientation of the edges of $H_k$, such that all
vertices of $H_k$ have out-degree~2.

It is clear that $x_0=1=y_0$ and Creed observes the following
recursion 
\[
\begin{array}{lcl}
x_k&=& 4x_{k-1}+2y_{k-1}\\
y_k&=& 4x_{k-1}+3y_{k-1}.\\
\end{array}
\]
This recursion formula can be verified by a simple
enumeration. In~\cite{creed07} a lemma from~\cite{v-tcocisrapg-01} is
used to argue 
that the sequence $x_k/y_k$ is non-repeating. In Lemma~\ref{lem:strimode} we
provide an easy argument to show from first principles that $x_k/y_k$ is
strictly monotonically decreasing.
 
Let $N_i$ denote the number of Eulerian orientations of $G_0$ which
have $i$ valid cross over boxes, that is $N_\ell$ is the number of Eulerian
orientations of $G$. The number of Eulerian orientations of $G_k$ is 
\[EO(G_k)=\sum_{i=0}^\ell N_ix_k^iy_k^{\ell-i}.\]
Hence the number $EO(G_k)/y_k^\ell$ is the value of the polynomial
$p(z)=\sum_{i=0}^\ell N_iz^i$ at the point $x_k/y_k$. Since computing
$y_k^\ell$ is easy the polynomial $p$ of degree 
$\ell$ can be evaluated at $\ell+1$ different points with $\ell+1$
calls to an oracle for counting Eulerian orientations of planar
graphs. Hence the coefficients of $p$ can be determined using polynomial
interpolation. In particular this yields a way to compute $N_\ell$. \qed 

\begin{lemma}\label{lem:strimode}
Let $x_0=1=y_0$ and $x_{k+1}=4x_k+2y_k$ while $y_{k+1}=4x_k+3y_k$. Then the
sequence $x_k/y_k$ is strictly monotonically decreasing.
\end{lemma}
\proof It follows directly from the recursion that
$x_{k+1}=y_{k+1}-y_k$ and $y_{k+1}=7y_k-4y_{k-1}$. We want to show
that
\[\frac{x_{k+1}}{y_{k+1}}=\frac{4x_{k}+2y_k}{4x_{k}+3y_k}<\frac{x_k}{y_k}
\Longleftrightarrow 4\left(\frac{x_k}{y_k}\right)^2-\frac{x_k}{y_k}-2>0\]
Since $x_k/y_k>0$ this inequality is satisfied if and only if
$x_k/y_k>(1+\sqrt{33})/8=:c$. Note that $1/(1-c)=3+4c$ can be easily
derived since $c$ solves $4t^2-t-2=0$.

It remains to show that $x_k/y_k=1-(y_{k-1}/y_k)>c$, which we do by
induction. Since $x_0/y_0=1$ and $x_1/y_1>6/7$ the induction base is trivial and the 
induction hypothesis guarantees that $1-c>y_{k-2}/y_{k-1}$ for $k\geq
2$. We have 

\[\frac{x_k}{y_k}=1-\frac{y_{k-1}}{y_k}>c\Longleftrightarrow 
\frac{y_k}{y_{k-1}}>\frac{1}{1-c}.\]
Using the induction hypothesis for the last inequality we obtain that
\[\frac{y_k}{y_{k-1}}
=\frac{7y_{k-1}-4y_{k-2}}{y_{k-1}}=
7-4\frac{y_{k-2}}{y_{k-1}}>3+4c=\frac{1}{1-c}\]\qed

The reduction from counting
perfect matchings in bipartite graphs in~\cite{mw-otnoeooag-96} to
counting Eulerian orientations creates
vertex degrees, which grow linearly with the number of vertices of the
graph. This implies also that the image of $\aa(v)=d(v)/2$ grows with
$n$. The following theorem shows that this dependence of $\aa$ on $n$ is not
necessary to obtain a  $\#{}P$-completeness result.

\begin{theorem}\label{thm:complexity}
For the following instances of graph classes and out-degree functions
$\aa$ the counting of $\aa$-orientations is $\#{}$P-complete.
\begin{enumerate}
\item Planar maps with $d(v)=4$ and $\aa(v)\in\{1,2,3\}$ for all
  $v\in V$. 
\item Planar maps with $d(v)\in\{3,4,5\}$ and $\aa(v)=2$ for all $v\in
  V$.
\end{enumerate}
\end{theorem}
The proof uses the planarization method from the proof of
Theorem~\ref{thm:creed} in conjunction with the 
following theorem  from~\cite{dl-atpogwlf-92}.
\begin{theorem}\label{thm:regmatch}
Counting perfect matchings of $k$-regular bipartite graphs is
$\#{}$P-complete for every $k\geq 3$.
\end{theorem}

\ni{\sl Proof of Theorem~\ref{thm:complexity}.} Perfect
matchings of a bipartite graph $G$ with vertex set $V=A\cup B$ are in
bijection with $\aa$-orientations of $G$ with $\aa(v)=1$ for $v\in A$
and $\aa(v)=d(v)-1$ for $v\in B$. This bijection is established  by
identifying matchings edges with edges directed from $A$ to $B$. Hence
in $k$-regular bipartite graphs counting perfect matchings  is equivalent
to counting what we call 1-$(k-1)$-orientations in the sequel.

We observe that the planarization method from the proof of
Theorem~\ref{thm:creed} can be used  in a more general setting.
Let $\mathcal{G}_D$ be the set of all graphs with vertex
degrees in $D\subset \NN$ and $\mathcal{P}_D$ the set of all planar
graphs with degrees in $D$. Let $I\subset\NN$ and associate with every
$G\in \mathcal{G}_D$ an out-degree function $\aa_G$ whose image is
contained in~$I$. Then, the proof of
Theorem~\ref{thm:creed} shows that counting the $\aa_G$-orientations
of graphs in $\mathcal{G}_D$ can be reduced to counting
$\aa_G'$-orientations of the graphs in $\mathcal{P}_{D\cup\{4\}}$
where the image of $\aa_G'$ is contained in $I\cup\{2\}$ for all
$G'\in\mathcal{P}_{D\cup\{4\}}$. 

When we apply this to 4-regular bipartite graphs with
1-3-orientations, that is perfect matchings, it yields the first claim
of the theorem since counting perfect matchings of bipartite 4-regular
graphs is $\#{}P$-complete by Theorem~\ref{thm:regmatch}. We give
two different proofs for the second claim. 

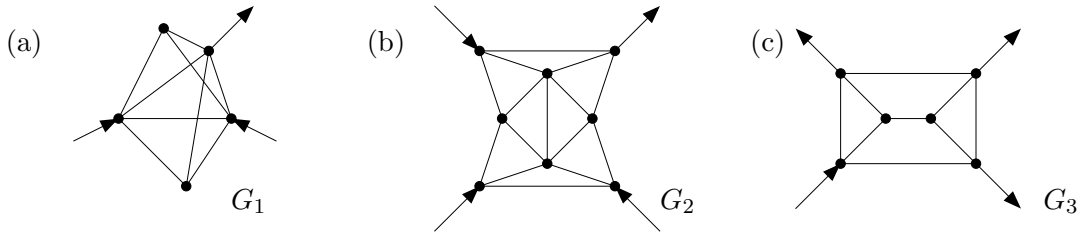
\PsFigCap{60}{alphagadgets}{Gadgets to translate (a) $\aa(v)=1$ to
  $\aa\equiv2$ for $d(v)=3$, (b) $\aa(v)=1$ to $\aa\equiv2$ for
  $d(v)=4$, and (c) $\aa(v)=3$ to $\aa\equiv2$ for $d(v)=4$. }
Let $G$ be a graph with a degree 3 vertex $v$ and $\aa_G$ an associated
out-degree function with $\aa_G(v)=1$. We substitute $v$ by the gadget
$G_1$ from Figure~\ref{fig:alphagadgets}~(a) to obtain a graph
$G'$. The gadget has five vertices that induce nine edges and 3 edges
connect it with the neighbors of $v$ in $G$. Let $\aa_{G'}=\aa_G$ on
$V(G)-v$ and be $\aa_{G'}(u)=2$ for $u\in V(G_1)$. Exactly one of the
connection edges must be directed away from $G_1$ in every
$\aa_{G'}$-orientation of $G'$. Note that $G_1$ is symmetric 
in the three connection vertices. It is easy to check that $G_1$ has
ten orientations with out-degree 2 at every vertex once the outgoing
connection edge has been chosen. Thus, every $\aa_G$-orientation is
associated with ten $\aa'_G$-orientations of $G'$ and since the
underlying $\aa_G$-orientations can be reconstructed from every
$\aa_{G'}$-orientation we obtain that $r_{\alpha_{G'}}(G')=10\cdot
r_{\alpha_G}(G)$.  

Let a 3-regular bipartite graph $G$ with vertex set $V=A\cup B$ be given and
$G'$ be obtained from $G$ by substituting every vertex $v\in A$ by a
copy of $G_1$. The number of 1-2-orientations of $G$ is
$10^{-|A|}\cdot r_{\alpha_{G'}}(G')$. Using the planarization method this
yields that counting 2-orientations of graphs from
$\mathcal{P}_{\{3,4,5\}}$ is  $\#{}P$-complete.

Similarly counting perfect matchings of 4-regular bipartite graphs can
be reduced to counting 2-orientations of graphs from
$\mathcal{P}_{\{3,4,5\}}$. We mention this second
proof since it uses planar gadgets. More precisely $G'$ is
obtained from $G$ by substituting the
vertices of one partition class of a 4-regular bipartite graph $G$
by the gadget $G_2$ shown
in  Figure~\ref{fig:alphagadgets}~(b) and the vertices from the other
partition class by the gadget $G_3$ shown in
Figure~\ref{fig:alphagadgets}~(c). Using a similar reasoning as above
one obtains that the number of 1-3-orientations of $G$ is $(26\cdot
6)^{-|A|}r_{\aa_{G'}}(G')$, where $G_2$ and $G_3$ give a blow
up factor $26$ respectively $6$. \qed

Having proved Theorem~\ref{thm:complexity} it is natural to ask whether
it is $\#{}P$-complete to count $\aa$-orientations for $k$-regular
planar graphs and constant $\aa$. This setting implies that $\aa\equiv
k/2$. Planar graphs have average degree less than 
six. Furthermore a 2-regular connected graph, that is a cycle, has two
Eulerian orientations. Hence the above question should be asked for
$k=4$. Note, that the planarization method yields vertices with $d(v)=4$
and $\aa(v)=2$, so we do not need to restrict our considerations to
planar graphs when trying to answer this question.\medskip

\ni\textbf{Problem.} Is it  $\#{}P$-complete to count Eulerian
orientations of 4-regular graphs?\medskip

\ni To the best of our knowledge even the following problem is also open.\medskip

\ni\textbf{Problem.} Is it $\#{}P$-complete to count Eulerian
orientations of graphs with degrees in $\{1,\ldots,k\}$, for some arbitrary but
fixed $k\in \NN$? \medskip

We present one more  $\#{}P$-completeness result since it has a nice
connection with the first problem stated above. Note, that for
4-regular bipartite graphs 2-factors and Eulerian orientations are in
bijection. For 3-regular bipartite graphs 2-factors are in bijection
with their complements the 1-factors, i.e. perfect matchings. Hence it
is $\#{}P$-complete to count 2-factors of 3-regular bipartite
graphs. The next theorem generalizes this observation.
\begin{theorem}\label{thm:2factor}  For every $i\geq 3,\; i\neq 4$
  counting 2-factors of $i$-regular bipartite graphs is $\#{}P$-complete. 
\end{theorem}
\proof The case $i=3$ follows from Theorem~\ref{thm:regmatch} as
explained above. The proof for $i\geq 5$ is a reduction from counting
2-factors of 3-regular bipartite graphs. The method comes from the
proof of Theorem~\ref{thm:regmatch} in~\cite{dl-atpogwlf-92}. 

The following preliminary considerations will be needed
later. We fix some edge $e_0$ of the complete bipartite graph on $2i$
vertices $K_{i,i}$. This graph has $i^2$ edges and every 2-factor
of $K_{i,i}$ has $2i$ edges. Let $c_i$ be the number of 2-factors of
$K_{i,i}$. We want to find the ratio between 2-factors of $K_{i,i}$
which contain $e_0$ and 2-factors which do not contain $e_0$. Consider
pairs of 2-factors and edges $(F,e)$. Obviously, there are $i^2\cdot c_i$
such pairs and we have that $|\{(F,e)\mid e\in
F\}|=2i\cdot c_i$ while $|\{(F,e)\mid e\not\in F\}|=(i^2-2i)\cdot
c_i$. It is obvious that 
\[
\{(F,e)\mid e\in F\}=\bigcup_{j=1}^{|E|}\{(F,e_j)\mid e_j\in F\}
\]
and symmetry implies, that all sets $\{(F,e_j)\mid e_j\in F\}$ have
the same cardinality. We conclude, that there are
$a_i=2c_i/i$ 2-factors including $e_0$. It follows similarly, that
there are $b_i=(1-2/i)c_i$ 2-factors not including $e_0$. We infer
that $a_i/b_i=2/(i-2)\leq 2/3$ for $i\geq5$. 

The gadget $P_i(k)$ is a concatenation of $k$ disjoint copies of
$K_{i,i}-e_0$, with $k-1$ connection edges as shown in
Figure~\ref{fig:bridgegadget} for $i=5$ and $k=4$. The gadget is
connected to the rest of the graph via two edges at the degree $i-1$
vertices.
%
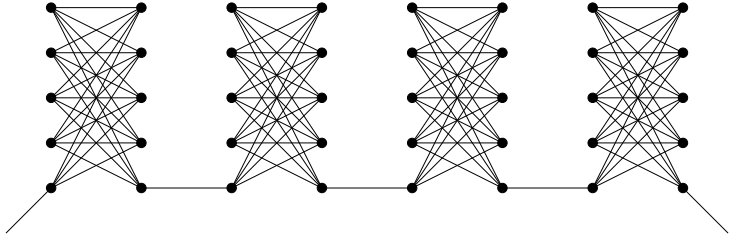
\PsFigCap{60}{bridgegadget}{The bridge gadget $P_5(4)$}

Let $G$ be a 3-regular bipartite graph with $2n$ vertices and $G'(k)$
be obtained from $G$ by augmenting it with $n(i-3)$ disjoint bridge
gadgets $P_i(k)$ such that $G'(k)$ is $i$-regular. Let $P$ be some fixed
bridge. Note that every 2-factor of $G'(k)$ includes all or none of the
connection edges of $P$. This is because a 2-factor is a partition
into cycles and therefore intersects every edge cut of cardinality two
in either none or both of the edges.

Let $c_G$ be the number of 2-factors of $G$. Among the 2-factors of
$G'(k)$ let $S$ denote the set of those, which induce a 2-factor of $G$,
and $S^c$ the other 2-factors. We have $|S|=c_G\cdot b_i^{k\cdot
  n(i-3)}$ since the 2-factors in $S$ cannot include any connection
edges and thus can be partitioned into a 2-factor of $G$ which is
augmented by one of the $b_i^k$ possible 2-factors in every bridge.

Next we want an upper bound for $|S^c|$. Every 2-factor in $S^c$
includes the connection edges of at least one bridge $P$. Since
$2^{3n}$ is the number of subsets of $E(G)$ each of which can be
augmented to a 2-factor in at most $a_i^k$ ways on $P$ we have that 
$|S^c|\leq a_i^kb_i^{k\cdot (n(i-3)-1)}2^{3n}$. The number of
2-factors of 
$G'(k)$ is $c_{G'(k)}=|S|+|S^c|$ for which we have the following
bounds. 

\[
\begin{array}{llllll}& c_G\cdot b_i^{k\cdot n(i-3)} &\leq& c_{G'(k)} &\leq&
  c_G\cdot b_i^{k\cdot 
  n(i-3)}+a_i^kb_i^{k\cdot (n(i-3)-1)}2^{3n}\\
 \Longleftrightarrow &c_G&\leq& c_{G'(k)} b_i^{-k\cdot n(i-3)}&\leq&
  c_G+2^{3n}\left(\frac{a_i}{b_i}\right)^k\end{array}\]   
Since $a_i/b_i\leq 2/3$ we have that 
\[2^{3n}\left(\frac{a_i}{b_i}\right)^k<\frac{1}{2}\mbox{ for }
k>3/(\log_23-1)n+1/(\log_23-1).\]
Note that the lower bound for $k$ is linear in $n$. 

We finally
conclude that $c_G= \lfloor c_{G'} b_i^{-k\cdot n(i-3)}\rfloor$ for 
$k$ large enough. Since  $b_i^{-k\cdot n(i-3)}$ is easy to compute and
$G_k$ has size polynomial in $n$ this proves the theorem. \qed

\ni\textbf{Remark.} We would like to point out, that the missing case
$k=4$ in Theorem~\ref{thm:2factor} cannot be fixed by substituting the
bridge gadget by another gadget. It is crucial that $a_i/b_i<1$ and
symmetry implies, that in
every 4-regular graph the number of 2-factors including a fixed edge
$e_0$ is equal to the number of 2-factors not containing $e_0$.\medskip 


\subsection{Approximation}
\label{subsec:perfmatch}


Let $f:V\rightarrow \NN$ be defined on the vertex set of a graph $G$. An
$f$-factor of a graph $G$ is a subgraph of $G$ such that every 
vertex has degree $f(v)$ in $F$.

Counting $\aa$-orientations can be reduced
to counting $f$-factors in bipartite planar graphs and to counting
perfect matchings in bipartite graphs.  We next describe these transformations.
They are useful because bipartite perfect
matchings have been the subject of extensive research (for
example~\cite{lovaszplummer,robertsonseymourthomas99,jsv-04}).

First, note that the $\aa$-orientations of $M$ are in bijection with
the the $\aa'$-orientations of the bipartite planar map $M'$ obtained from
$M$ by subdividing every edge once. Here $\aa'(v)=\aa(v)$ for the
original vertices of $M$ and $\aa'(v)=1$ for all subdivision
vertices. The $\aa'$-orientations of $M'$ are in bijection with the
$f$-factors of $M'$ where $f(v)=\aa'(v)$ for all vertices of $M'$. The
bijection works by identifying factor edges with edges directed from a
vertex of $M$ to an edge vertex.

The idea for the next transformation is due to
Tutte~\cite{t-aspotftffg-54}.
The graph $M'$ gets blown up into a graph $M''$ such that $M''$ has
$\prod_{v\in M}(d(v)-f(v))!$ times as many perfect matchings, i.e.,
1-factors, as there are $f$-factors of $M'$.  To obtain $M''$ from
$M'$ substitute $v\in V(M)$ by a $K_{d(v),d(v)-f(v)}$, such that each
of the $d(v)$ edges incident to $v$ in $M'$ connects to one of the
vertices from the partition class of cardinality $d(v)$. 

In~\cite{jsv-04} Jerrum, Sinclair, and Vigoda give a
fully polynomial randomized approximation scheme for counting perfect
matchings of bipartite graphs. Thus, the above transformation yields a
fully polynomial randomized approximation scheme for
$\alpha$-orientations as well.  

The number of perfect matchings of a bipartite graph with
a Pfaffian orientation can be computed in polynomial
time. Little~\cite{little} gave a full characterization of graphs with
a Pfaffian orientation and in~\cite{robertsonseymourthomas99} a
polynomial time algorithm to test whether a given graph is Pfaffian is
introduced. As a special case of Little's characterizations it follows that
all planar graphs are  
Pfaffian. Hence, in all cases where $M''$ is planar the counting is
easy. For spanning trees the above transformations yields planar
graphs, while for Eulerian orientations it does not, as Theorem~\ref{thm:creed}
implies. Although we do not have a hardness result for Schnyder woods
and bipolar orientations on planar maps, there are in both cases
instances for which the transformation yields a
non-Pfaffian graph. One such example is the augmented triangular grid
from Section~\ref{sec:swont}. Figure~\ref{fig:bipolar-nopfaffian}
shows a local structure with five vertices and four faces that implies
that $M''$ is not Pfaffian. The figure shows simplified versions of
$\widehat{M}'$ and $\widehat{M}''$. We can choose
$\widehat{M}'=\widehat{M}$ since the angle graph is bipartite. The
Tutte transformation substitutes face vertices by $K_{1,3}$ and
primal vertices of degree $d$ by $K_{2,d}$. Instead one can simply create a
copy of every primal vertex with the same neighborhood as the original
and leave the face vertices unchanged to obtain a simplified version
of $\widehat{M}''$.

%
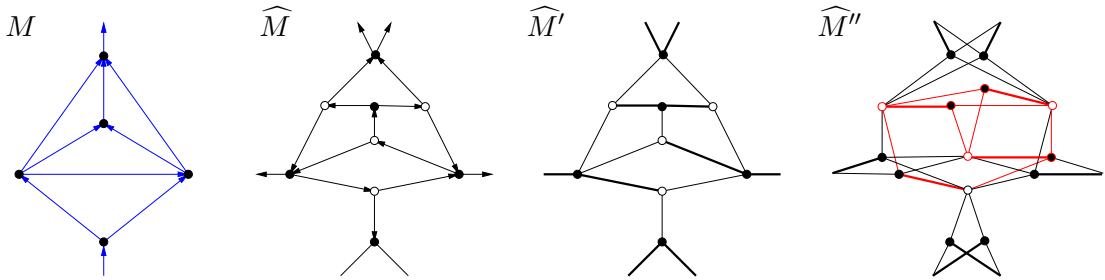
\PsFigCap{45}{bipolar-nopfaffian}{An oriented subgraph of a
triangulation with a bipolar orientation that induces a central
$K_{3,3}$-subdivision.} 

\section{Conclusions}
In this paper we have studied the maximum number of
$\alpha$-orientations for different classes of planar maps and
different $\alpha$. In most cases we have exponential upper and lower bounds
$c_L^n$ and $c_U^n$ for this number. 

The obvious problem is to improve on the constants
$c_L$ and $c_U$ for the different instances. We think, that in
particular improving the upper bound of $8^n$ for the number of
Schnyder woods on 3-connected planar maps is worth further
efforts. For bipolar orientations the more efficient encoding from 
Theorem~\ref{thm:plusminus} helps to improve the upper bound. We think
that finding a more efficient encoding for Schnyder woods  might be
needed to substantially improve on the $8^n$ bound.

Results by Lieb~\cite{l-treosi-67} and
Baxter~\cite{b-fmtl-69} yield the exact asymptotic behavior of the
number of Eulerian orientations for the square and triangular grid on
the torus. This yields upper bounds for the number of 2-orientations
on the square grid and the number of Schnyder woods on the triangular
grid. We have not been able to take
advantage of these results for improving the lower bounds for the
number of 2-orientations respectively Schnyder woods. 

For some instance of $\aa$-orientations there are $\#{}P$-completeness
results. This contrasts with spanning trees
and planar bipartite perfect matchings for which polynomial algorithms
are available. It remains open to determine the complexity of counting
Schnyder woods and bipolar orientations on planar maps and of counting
Eulerian orientations of graphs with bounded maximum degree.\medskip

\ni\textbf{Acknowledgments.} We would like to thank Graham Brightwell for
interesting discussions and valuable hints in connection with Lieb's
3-coloring of  
the square grid. We thank Christian Krattenthaler for directing us to 
reference~\cite{kpw-tam-00} and Mark Jerrum for bringing P\'aid\'i
Creed's work to our attention. Thanks to Pa\'id\'i Creed for
sending us a preliminary version of his proof. Florian Zickfeld was
supported by the Studienstiftung des deutschen Volkes.


\bibliographystyle{my-siam}
\bibliography{noor-rev}

\end{document}